\documentclass{article}
\date{14.4.2012}  
\usepackage{hyperref}
\hypersetup{
    urlcolor = red,
	colorlinks = true,
	linkcolor = blue,
	citecolor = blue,
    linktocpage = true,
	pdftitle = {Making Lifting Obstructions Explicit}, 
    pdfauthor = {Karl-Hermann Neeb, Friedrich Wagemann and Christoph Wockel},
    bookmarksopen = true,
    bookmarksopenlevel = 1,
    hypertexnames = false
  }

\usepackage[american]{babel}
\usepackage{amsfonts,amsmath,amssymb}

\usepackage{amsthm}

\usepackage[all,cmtip,line]{xy}
\xyoption{all}
\input{xypic}
\newdir{ >}{{}*!/-8pt/\dir{>}}
\newcommand{\susp}[1]{\ensuremath{\mathop{{\rm sp}_{#1}}}}

\renewcommand{\uline}{\underline}
\newcommand{\oline}{\overline}
\newcommand{\into}{\hookrightarrow} 
\newcommand{\onto}{\to\mskip-14mu\to}
\newcommand{\cU}{\mathcal U}

\renewcommand{\cH}{\mathcal H}
\newcommand{\fz}{\mathfrak z}
\newcommand{\fk}{\mathfrak k}

\newcommand{\subeq}{\subseteq} 
\newcommand{\eset}{\emptyset} 
\newcommand{\id}{{\rm id}} 
\renewcommand{\:}{\colon} 
\newcommand{\res}{\mid} 
\renewcommand{\hat}{\widehat} 
\newcommand{\pr}{\mathop{{\rm pr}}\nolimits}
\newcommand{\inn}{\mathop{{\rm int}}\nolimits}
\newcommand{\Ext}{\mathop{{\rm Ext}}\nolimits}
\newcommand{\Tor}{\mathop{{\rm Tor}}\nolimits}
\newcommand{\End}{\mathop{{\rm End}}\nolimits}
\renewcommand{\tilde}{\widetilde}
\newcommand{\trile}{\trianglelefteq}
\newcommand{\Ch}{\mathop{{\rm Ch}}\nolimits}
\newcommand{\obs}{\mathop{{\rm obs}}\nolimits}
\newcommand{\sing}{\mathop{{\rm sing}}\nolimits}

\newcommand{\Homeo}{\mathop{{\rm Homeo}}\nolimits}

\newcommand{\PU}{\mathop{{\rm PU}}\nolimits}
\newcommand{\SU}{\mathop{{\rm SU}}\nolimits}
\newcommand{\U}{\mathop{{\rm U}}\nolimits}

\newcommand{\im}{\mathop{{\rm im}}\nolimits}
\newcommand{\Hom}{\mathop{{\rm Hom}}\nolimits}
\newcommand{\Bun}{\mathop{{\rm Bun}}\nolimits}
\newcommand{\coker}{\mathop{{\rm coker}}\nolimits}
\renewcommand{\1}{\mathbf 1}

\newcommand{\ssssarr}{\hbox to 15pt{\rightarrowfill}}
\newcommand{\sssarr}{\hbox to 20pt{\rightarrowfill}}
\newcommand{\ssarr}{\hbox to 30pt{\rightarrowfill}}
\newcommand{\sarr}{\hbox to 40pt{\rightarrowfill}}
\newcommand{\arr}{\hbox to 60pt{\rightarrowfill}}
\newcommand{\larr}{\hbox to 60pt{\leftarrowfill}}
\newcommand{\Arr}{\hbox to 80pt{\rightarrowfill}}

\newcommand{\sssmapright}[1]{\smash{\mathop{\sssarr}\limits^{#1}}}
\newcommand{\ssmapright}[1]{\smash{\mathop{\ssarr}\limits^{#1}}}
\newcommand{\smapright}[1]{\smash{\mathop{\sarr}\limits^{#1}}}

\def\bS{\mathbb{S}}
\def\bB{\mathbb{B}}
\def\R{\mathbb{R}}
\def\N{\mathbb{N}}
\def\C{\mathbb{C}}

\def\Z{\mathbb{Z}}

\def\T{\mathbb{T}}

\newtheorem{theo}{Theorem}[section]
\newtheorem{defi}[theo]{Definition}
\newtheorem{exa}[theo]{Example}
\newtheorem{rema}[theo]{Remark}
\newtheorem{prop}[theo]{Proposition}
\newtheorem{cor}[theo]{Corollary}
\newtheorem{lem}[theo]{Lemma}
\newtheorem{pro}[theo]{Problem} 

\newenvironment{defn}{\begin{defi}\rm}{\end{defi}}
\newenvironment{rem}{\begin{rema}\rm}{\end{rema}}
\newenvironment{ex}{\begin{exa}\rm}{\end{exa}}
\newenvironment{prob}{\begin{pro}\rm}{\end{pro}}

\newenvironment{prf}{\begin{proof}}{\end{proof}}
\newcommand{\mlabel}[1]{\marginpar{#1}\label{#1}}
\renewcommand{\mlabel}{\label}

\begin{document}
\title{Making Lifting Obstructions Explicit} 

\author{Karl-Hermann Neeb\\
        Universit\"at Erlangen-N\"urnberg\\
\and Friedrich Wagemann\\
     Universit\'e de Nantes\\
\and Christoph Wockel\\
     Universit\"at Hamburg}

\maketitle

\begin{abstract}
If $P \to X$ is a topological principal $K$-bundle 
and $\hat K$ a central extension of $K$ by $Z$, 
then there is a natural obstruction class $\delta_1(P) \in \check H^2(X,\uline Z)$ 
in sheaf cohomology whose vanishing is equivalent to the existence 
of a $\hat K$-bundle $\hat P$ over $X$ with $P \cong \hat P/Z$. 
In this paper we establish a link between homotopy theoretic 
data and the obstruction class $\delta_1(P)$ which in many cases can 
be used to calculate this class in explicit terms. 
Writing $\partial_d^P \: \pi_d(X) \to \pi_{d-1}(K)$ for the 
connecting maps in the long exact homotopy sequence, two of our 
main results can be formulated as follows. 
If $Z$ is a quotient of a contractible group by the discrete group $\Gamma$, 
then the homomorphism $\pi_3(X) \to \Gamma$ induced by 
$\delta_1(P) \in \check H^2(X,\uline Z) \cong H^3_{\rm sing}(X,\Gamma)$ 
coincides with $\partial_2^{\hat K} \circ \partial_3^P$ 
and if $Z$ is discrete, then 
$\delta_1(P) \in \check H^2(X,\uline Z)$ 
induces the homomorphism $-\partial_1^{\hat K} \circ \partial_2^P \: \pi_2(X) \to Z$. 
We also obtain some information on obstruction classes defining trivial 
homomorphisms on homotopy groups. 
\end{abstract}

\tableofcontents

\section{Introduction}

Let $K$ be a topological group and $q \: P \to X$ be a $K$-principal bundle 
on the paracompact Hausdorff space~$X$. 
Suppose further, that $\hat K$ is a topological  central extension 
of $K$ by $Z$ (which is always assumed to be a 
locally trivial $Z$-bundle over $K$). 
We say that {\it $P$ lifts to a $\hat K$-bundle} if there exists a 
$\hat K$-principal bundle $\hat P$ over $X$ with 
$\hat P/Z \sim P$ (as $K$-bundles). 
By local trivializations and the corresponding transition functions, 
we obtain a bijection between the set $\Bun(X,K)$ of equivalence 
classes of topological principal $K$-bundles over the paracompact space $X$ and 
the sheaf cohomology set 
$\check H^1(X,\uline K)$, where 
$\uline K$ is the sheaf of germs of continuous $K$-valued 
functions on $X$. Now the short exact sequence 
$Z \into \hat K \to K$ of topological groups 
induces a short exact sequence of sheaves of groups 
\[ \1 \to \uline Z \to \uline{\hat K} \to \uline K \to \1,\] 
which in turn leads to an exact sequence in sheaf cohomology 
\begin{eqnarray}\label{eq:shexact}
& C(X,K) \cong \check H^0(X,K) \ssmapright{\delta_0} 
\Bun(X,Z) \cong \check H^1(X,\uline Z)\\
&\to \Bun(X,\hat K) \cong \check H^1(X,\uline{\hat K}) 
\to \Bun(X,K)\cong \check H^1(X,\uline K) 
\ssmapright{\delta_1} \check H^2(X,\uline Z)\notag
\end{eqnarray}
(cf.\ \cite{Gr58}). 
Note that for non-abelian $K$ and $\hat{K}$ the latter is just an exact sequence of
pointed sets rather than abelian groups.
For the cohomology class 
$[P] \in \check H^1(X,\uline{K})$ representing the 
$K$-bundle $P$, the 
class $\delta_1([P]) \in \check H^2(X,\uline Z)$ obtained from 
the connecting map can therefore be interpreted as 
an {\it obstruction class}. It vanishes if and only 
if $P$ lifts to a $\hat K$-bundle. 
In the present paper we address the problem of making this obstruction 
class as explicit as possible in terms of 
characteristic data attached to the central extension $\hat K$ and the bundle~$P$.

To this end, we assume that the identity 
component $Z_0$ of $Z$ is open and a {\it $K(\Gamma,1)$-group}, i.e., 
locally contractible, $\pi_1(Z)$ is isomorphic to the discrete abelian group 
$\Gamma$, and its simply 
connected covering group $\tilde Z_0$ is contractible 
and divisible (cf.\ Definition~\ref{def:kgrp}).
\begin{footnote}{Typical examples of $K(\Gamma,1)$-groups are quotients 
of a topological vector space $\fz$ by a discrete subgroup $\Gamma$.}  
\end{footnote}
This implies in particular that $Z_0$ is divisible, so that 
$Z$ is a direct product $Z \cong Z_0 \times D$, where 
$D \cong \pi_0(Z)$ is a discrete abelian group. 
From the short exact sequence $\Gamma \to \tilde Z_0 \to Z_0$ 
and the corresponding short exact sequence 
\[ 0 \to \uline \Gamma \to \uline{\tilde Z_0} \to \uline {Z_0} \to 0 \] 
of abelian group sheaves we obtain the corresponding long exact sequence 
in sheaf cohomology. For each paracompact space $X$, 
the vanishing of the cohomology 
$\check H^n(X,\uline{\tilde Z_0})$ for $n > 0$ 
(\cite[Prop.~4]{Hu61}) 
therefore yields natural isomorphisms 
\begin{equation}
  \label{eq:iota}
\iota_n \: \check H^n(X,\uline {Z_0})\to 
\check H^{n+1}(X,\uline \Gamma) \cong \check H^{n+1}(X,\Gamma), 
\end{equation}
where the latter isomorphism follows from the discreteness of 
$\Gamma$ (\cite[Prop.~4]{Hu61}, resp.\ \cite[Ex.~5.1.10, p.~230]{God}). 
Therefore our problem boils down to describe a cohomology class in 
\[ \check H^2(X,\uline Z) 
\cong \check H^2(X,\uline{Z_0}) \oplus \check H^2(X,\uline{D}) 
\cong \check H^3(X,\Gamma) \oplus \check H^2(X,D). \] 
In particular, we can study the case of a connected group $Z = Z_0$ and a discrete 
group $Z = D$ independently. 

Throughout we shall assume that $X$ is paracompact (needed for the homotopy theory 
of bundles) and that $X$ is locally contractible, which leads to a natural isomorphism 
$\check H^k(X,\Gamma) \to H^k_{\rm sing}(X,\Gamma)$ between \v Cech cohomology and singular 
cohomology. Since we can evaluate singular cohomology classes on homology classes, 
we obtain for each  $n \in \N_0$ homomorphisms 
\begin{equation}
  \label{eq:alphan0}
\alpha_n \: \check H^n(X,\Gamma) \to \Hom(\pi_n(X),\Gamma) 
\end{equation}
of abelian groups (see Definition \ref{defalpha} for details). 
Therefore our first step in analyzing the obstruction classes consists in 
identifying their image under $\alpha_2$, resp., $\alpha_3$. 
To formulate our results, we write 
\begin{equation}
  \label{eq:conmap}
\partial^P_j \:  \pi_j(X) \to \pi_{j-1}(K) 
\end{equation}
for the connecting map in the long exact homotopy sequence of the 
bundle $P \to X$. 
Then our main result can be stated as follows: 

\begin{theo} \mlabel{thm:main} 
Let $X$ be a paracompact locally contractible Hausdorff space, 
$K$ a connected locally contractible topological group, $P \to X$ a principal $K$-bundle, and 
$\hat K$ be a central extension of $K$ by $Z$. Then the following assertions 
hold: 
\begin{description}
\item[\rm(a)] If $Z$ is a $K(\Gamma,1)$-group, then $\alpha_3(\iota_2(\delta_1([P]))) 
= \partial_2^{\hat K} \circ \partial_3^P$. 
\item[\rm(b)] If $Z = D$ is discrete, then $\alpha_2(\delta_1([P])) 
= -\partial_1^{\hat K} \circ \partial_2^P$. 
\end{description}
\end{theo}

Since $\pi_2(K)$ vanishes for every finite dimensional 
Lie group $K$ (\cite{Ca36}), part (a) of the preceding theorem 
shows that for these groups the obstruction class defines the trivial 
homomorphism $\pi_3(X) \to \pi_1(Z)$. 
If $X$ is $2$-connected, then $H^3(X,\Gamma) \cong \Hom(\pi_3(X),\Gamma)$, 
so that (a) completely determines the obstruction class. 
If $X$ is $1$-connected, then $H^2(X,\Gamma) \cong \Hom(\pi_2(X),\Gamma)$ 
and (b) likewise determines the obstruction class for $Z = D$. 

For the proof of these results we proceed as follows. 
We start in Section~\ref{sec:0} 
with some generalities on obstruction classes. 
Since we want to express these obstruction classes in terms of easily 
accessible data attached to the central extension $\hat K$ and the bundle $P$, 
it is important to understand which kind of data is 
relevant for our purpose. 
Since we are mainly interested in the case where $K$ is a Lie group and 
$Z_0$ is a quotient of a topological vector space, we show in 
Section~\ref{sec:liegrp} that, for a connected locally contractible 
topological group $K$,  
the relevant data related to $\hat K$ consists of the 
homomorphisms 
\[ \partial_1^{\hat K} \: \pi_1(K) \to \pi_0(Z) \quad \mbox{ and } \quad 
\partial_2^{\hat K} \: \pi_2(K) \to \pi_1(Z), \] 
and, if $\partial_2^{\hat K}$ vanishes, the abelian extension 
\[ [\pi_1(\hat K)] \in \Ext(\pi_1(K), \pi_1(Z)). \]  
The obstruction class $\delta_1([P])$ depends only on this data. 
In the context of Lie group extensions, $\partial_2^{\hat K}$ 
can be written as the period homomorphism of a left invariant 
$2$-form (\cite[Prop.~5.11]{Ne02}) which is a crucial key for its 
explicit determination and hence for the calculation of obstruction 
classes. 

In Section~\ref{sec:sn} we collect 
various results on spheres. In particular, we describe an 
explicit way to calculate \v Cech cohomology classes 
with values in discrete abelian groups, introduce the 
suspension homomorphism and recall how $K$-bundles on 
$\bS^n$ are classified by elements of $\pi_{n-1}(K)$. 
Then we turn in Section~\ref{sec:1} to the 
case where $Z = D$ is discrete and prove 
Theorem~\ref{thm:main}(b). 
As a byproduct, we obtain an analog of 
Theorem~\ref{thm:main} for 
$\delta_0 \: C(X,K) \to \check H^1(X,\uline Z)$ 
(Proposition~\ref{prop:lowdeg}). 

In Section~\ref{sec:redux} we reduce the proof of 
Theorem~\ref{thm:main}(a) first to the case 
$X=\bS^3$ by pulling back geometric data by continuous maps 
$\bS^3 \to X$ and then we complete the proof by verifying the 
theorem for $X = \bS^3$ by direct computation. 
In the course of the proof, the following diagram plays a central 
role: 

\begin{equation}
  \label{eq:kdiag}
\vspace{.5cm}
\hspace{0cm}
\xymatrix{
\check{H}^0(\bS^2,\underline{K}) \cong C(\bS^2, K)
 \ar[r]^{\quad\quad\delta_0}\ar[d]^{} 
&  \check{H}^1(\bS^2,\underline{Z}) 
\ar[d]^{-s} \\
\check{H}^1(\bS^3,\underline{K})\cong \Bun(\bS^3,K) \ar[r]_{\quad\quad
\delta_1} & 
\check{H}^2(\bS^3,\underline{Z})}
\end{equation}

\noindent Here the vertical maps are suspension maps, 
so that the commutativity of the diagram expresses a compatibility of 
connecting maps with suspensions. 
  
In Section~\ref{sec:6}, we discuss situations where 
the obstruction classes lie in the subgroup 
$\Lambda^3(X,\Gamma) := \ker \alpha_3 \subeq \check H^3(X,\Gamma)$ 
of {\it aspherical cohomology classes} 
inducing the zero homomorphism $\pi_3(X) \to \Gamma$. 
If $X$ is $1$-connected, then 
$\Lambda^3(X,\Gamma) \cong \Ext(\pi_2(X),\Gamma)$. 
Here our main result is Theorem~\ref{thm:7.11} which provides 
complete information on the obstruction class 
for $\partial_2^P = 0$ and a flat extension of $K$ defined by a 
homomorphism $\gamma \: \pi_1(K) \to Z$, which implies that 
$\delta_1([P]) \in \Ext(H_2(X),\Gamma) \subeq \Lambda^3(X,\Gamma)$.
We also explain how to construct other types of aspherical obstruction 
classes.

\noindent{\bf Acknowledgment:}
FW thanks TU Darmstadt and Universit\"{a}t Hamburg for 
supporting research visits of FW during 
which parts of this paper have been developed. 
In the same way, CW thanks the program ``Mathematiques en Pays de la Loire''
for sponsoring a research visit in Nantes.
We thank Thomas Nikolaus for informing us about the suspension
morphism in \v{C}ech cohomology and Linus Kramer for pointing us to the
reference \cite{Cal09} in Remark \ref{rem:5.6}.     

\section*{Conventions}

Unless stated otherwise
$X$ and $Y$ are paracompact locally
contractible Hausdorff spaces, 
$\Gamma$ is a discrete abelian group and
$Z$ is a locally contractible abelian topological group whose 
identity component $Z_0$ is a $K(\Gamma,1)$-group. 
Furthermore, $\mathfrak{z}$ always is a topological
vector space and when writing $\mathfrak{z}/\Gamma$ we always assume that
$\Gamma\subset \mathfrak{z}$ is a discrete \emph{subgroup}. 
We assume that we
are given a central extension
\begin{equation}\label{ext1}
 Z\into \hat{K}\to K,
\end{equation}
where $K$ is a locally contractible topological group and 
\eqref{ext1} is required
to be a locally trivial principal bundle. The set of equivalence classes 
of such extensions form a group $K$ is denoted by $\Ext_c(K,Z)$. 
If, in addition, $K$ is a Lie group we write $\Ext_{s}(K,Z)$ for the group 
of smooth central extensions, i.e., central extensions for which the group $\hat K$ in 
 \eqref{ext1} is a Lie group and a smooth principal $Z$-bundle. 
If $K$ and $Z$ are discrete, then $H^{n}_{{\rm grp}}(K,Z)\cong H^{n}(BK,Z)$
denotes the ordinary group cohomology. We will also sometimes use the
locally continuous cohomology groups $H^{n}_{c}({K,Z})$ and the locally smooth
cohomology groups $H^{n}_{s}(K,Z)$ in case that $K$ is a Lie group (cf.\ \cite{WW}). 
Whilst we are using different kinds of cohomology groups
we will throughout only use singular homology of spaces, which we denote by
$H_{n}(X)$. We also recall that the \v Cech cohomology groups 
$\check H^n(X,\Gamma)$ and the singular cohomology groups 
$H^n(X,\Gamma) := H^n_{\sing}(X,\Gamma)$ 
are isomorphic because $X$ was 
assumed to be locally contractible (\cite[Thm.~5.10.1, p. 228 and Ex.~3.9.1, 
p.~159]{God}). 

We recall that, for $k =0,1$ (and for $k \in \N$ if $\hat K$ is abelian), 
the connecting homomorphisms 
\begin{equation*}
 \delta_{k}\: \check{H}^{k}(X,\uline{K})\to \check{H}^{k+1}(X,\uline{Z})
\end{equation*}
can be constructed by representing a cohomology class $c$ by a cocycle 
$(g_{i_{0},\ldots, i_{k}})$ on an open cover $\cU = (U_i)_{i \in I}$ of $X$ for which 
all functions $g_{i_{0},\ldots, i_{k}}$ have continuous lifts 
leading to a $\hat{K}$-valued cochain $\hat g_{i_{0},\ldots,i_{k}}$. Then  $\delta_k(c)$ is the cohomology class represented by the $\uline Z$-valued cocycle 
$\delta (\hat{g})_{i_{0},\ldots,i_{k+1}}$.

\section{Generalities on the obstruction class} \mlabel{sec:0}

In this section we start with some general remarks on 
the obstruction class $\delta_1([P])$. We write  
\[ \obs^{\hat K}([P]) := \obs_P([\hat K]) := \delta_1([P]) 
\in \check H^2(X,\uline Z) \] 
for the corresponding obstruction class and recall 
that the group structure on $\Ext_c(K,Z)$ is defined by the Baer sum
\[ [\hat K_1] \oplus [\hat K_2] := [\hat K_3] 
\quad \mbox{ with } \quad 
\hat K_3 := (\hat K_1 \times_K \hat K_2)/\{ (z,z^{-1}) \: z \in Z \}.\]  

\begin{lem} \mlabel{lem:0.1} {\rm(a)} The map 
\[ \obs_P \: \Ext_c(K,Z) \to \check H^2(X,\uline Z) \] 
is a group homomorphism. 

{\rm (b)} If $K$ is discrete, then $\obs_P$ factors through a homomorphism 
\[ \obs_P^d \: \Ext_c(K,Z) \to \check H^2(X,Z_d),\] 
where $Z_d$ denotes $Z$, endowed with the discrete topology. 
\end{lem}

\begin{prf} (a) Let $q_j \: \hat K_j \to K$ be two 
$Z$-extensions of $K$, 
$(g_{k\ell}) \in \check Z^1(X,\uline K)$ be 
a $1$-cocycle and 
$(h_{k\ell}) \in \check C^1(X,\uline{\hat K_1})$,  
$(h'_{k\ell}) \in \check C^1(X,\uline{\hat K_2})$ be lifts of 
$(g_{k\ell})$. 
Then $h''_{k\ell} := (h_{k\ell}, h'_{k\ell})$ 
is a lift of $(g_{k\ell})$ with values in the 
$Z \times Z$-extension of $K$ given by the fiber product 
$\hat K_1 \times_K \hat K_2$ and 
$\delta h'' = (\delta h, \delta h') \in \check Z^2(X,\uline{Z^2})$ 
represents $(\obs_P([K_1]), \obs_P([K_2]))$. 

The central extension $[\hat K_1] \oplus [\hat K_2]\in \Ext_c(K,Z)$ 
is defined as the quotient $\hat K_3$ of $\hat K_1 \times_K \hat K_2$ 
by the central subgroup $\{ (z,z^{-1}) \: z \in Z \}$, 
which is the kernel of the multiplication homomorphism 
$\mu_Z \: Z^2 \to Z$. From that it follows that 
$\obs_P([\hat K_3]) = \obs_P([\hat K_1]) 
+ \obs_P([\hat K_2]).$

(b) This follows immediately from the bijection
$\Ext_c(K,Z_d) \to \Ext_c(K,Z)$ induced from the continuous bijection 
$Z_d \to Z$. 
\end{prf}

From the naturality of the connecting maps in 
\eqref{eq:shexact}, we immediately obtain the naturality 
of the obstruction classes: 
 
\begin{lem} \mlabel{lem:0.2} 
{\rm(a)} If $f \: Y \to X$ is a continuous map, $P \to X$ a $K$-bundle and 
$\hat K$ a $Z$-extension of $K$, then 
\[ \delta_1([f^*P]) = f^*\delta_1([P]) \in 
\check H^2(Y, \uline Z). \] 

{\rm(b)} If $\phi \: K_1 \to K_2$ is a continuous morphism of topological groups, 
\[\phi^* \: \Ext_c(K_2,Z) \to \Ext_c(K_1,Z), \quad [\hat K_2] \mapsto [\phi^*\hat K_2],\] 
$P \to X$ a $K_1$-bundle and $\phi_*P := P \times_\phi K_2$ the associated $K_2$-bundle, 
then 
\[ \obs_P([\phi^*\hat K_2]) = \obs_{\phi_*P}([\hat K_2]).\] 

{\rm(c)} If $\phi \: Z_1 \to Z_2$ is a continuous morphism of abelian groups 
and 
\[\phi_* \: \Ext_c(K,Z_1) \to \Ext_c(K,Z_2), \quad [\hat K] \mapsto [\phi_*\hat K],
\quad \phi_*\hat K = \hat K \times_\phi Z_2,\] 
then 
\[ \phi_* \obs_P([\hat K]) = \obs_{P}([\phi_*\hat K]).\] 
\end{lem}

\begin{defn} For a principal $K$-bundle $q_P \: P \to X$ over the topological space 
$X$, we write 
\[ \partial^P_n \: \pi_n(X) \to \pi_{n-1}(K) \] 
for the {\it connecting homomorphisms} in the long exact homotopy sequence 
of~$P$. Since we shall need this information below, we recall how these
maps are constructed. First we choose base points 
$p_0 \in P$ and $x_0 \in X$ with $q_P(p_0) = x_0$ and observe that 
the map $\iota \: K \to P_{x_0}, k \mapsto p_0.k$ is a homeomorphism. 
Let $I := [0,1]$ be the unit interval. 
For a continuous based map $\sigma \: \bS^n \cong I^n/\partial I^n \to X$ 
we choose a continuous lift $\tilde\sigma \: I^n \to P$. 
Then $\tilde\sigma(\partial I^n) \subeq P_{x_0}$, and we put 
\[ \partial_n^{P}([\sigma]) 
= [\iota^{-1} \circ \tilde\sigma\mid_{\partial I^n}] \in \pi_{n-1}(K).\] 
\end{defn}

As explained in the following remark, the connecting map 
$\partial_1^P$ is the obstruction for the reduction of the structure group 
to the identity component $K_0$. 
As Theorem~\ref{thm:main} shows, the higher 
connecting maps $\partial_2^P$, resp., $\partial_3^P$ are closely related to the 
obstruction classes for discrete, resp., connected groups~$Z$.

\begin{rem} (Geometric interpretation of $\partial_1^P$) 
Let $K_0 \subeq K$ be the identity component. Then 
$BK_0 := EK/K_0$ is a classifying space for $K_0$. As $K_0$ is connected, 
$BK_0$ is simply connected. Moreover, the natural map 
$BK_0 \to BK \cong EK/K$ is a covering because $K_0$ is open in $K$. 
Therefore $BK_0$ is the universal covering space of $BK$. 

Since a continuous map $f \: X \to BK$ can be lifted to $\tilde{BK}$ if and only if 
$\pi_1(f) = \partial_1^{f^*EK} \: \pi_1(X) \to \pi_1(BK) \cong \pi_0(K)$ 
vanishes, it follows that a $K$-bundle $P \to X$ 
has a reduction to a $K_0$-bundle if and only if the connecting homomorphism 
$\partial_1^P \: \pi_1(X) \to \pi_0(K)$ is trivial. 
\end{rem}

In Remark~\ref{rem:1conn}(c) below we shall see that the vanishing of the connecting 
map $\partial_2^P \: \pi_2(X) \to \pi_1(K)$ also has a simple geometric interpretation 
because it is equivalent to the existence of a $\tilde K$-lift of the pullback of the 
$K$-bundle $P$ to the universal covering space of~$X$. 
A geometric interpretation of the vanishing of $\partial_3^P$ 
can be found in Remark~\ref{rem:7.21}(b). 

\begin{lem} \mlabel{lem:0} If $P \to X$ is a $K$-bundle and 
$f \: (Y,y_0) \to (X,x_0)$ is a continuous based map, then 
\[ \partial^P_n \circ \pi_n(f,y_0) = \partial^{f^*P}_n \: \pi_n(Y,y_0) 
\to \pi_{n-1}(K).\] 
\end{lem}

\begin{prf} Pick a base point $p_0 \in P$ over $x_0$. 
For a continuous based map \break 
$\sigma \: \bS^n \cong I^n/\partial I^n \to Y$ 
we have 
\[ \partial_n^{f^*P}([\sigma]) = [\iota^{-1} \circ 
\tilde\sigma\mid_{\partial I^n}],\] 
where $\tilde\sigma \: I^n \to f^*P$ is a continuous base point preserving 
lift of $\sigma$. 
Since the projection $\pr_P \: f^*P \to P$ preserves base points, 
$\pr_P \circ \tilde\sigma$ is a base point preserving lift of 
$f \circ \sigma$, so that 
\[ \partial_n^P([f\circ \sigma]) 
= [\iota^{-1} \circ \pr_P \circ \tilde\sigma\mid_{\partial I^n}] 
= \partial^{f^*P}([\sigma]).\] 
This proves our assertion. 
\end{prf}

\begin{rem} \mlabel{rem:0.2} 
Since $Z_0$ is open in $Z$ and divisible, 
$Z \cong Z_0 \times D$ is a direct product of the discrete group 
$D := Z/Z_0 =: \pi_0(Z)$ and the connected group $Z_0$. 
This product decomposition leads to 
\[ \Ext_c(K,Z) \cong \Ext_c(K,Z_0) \oplus \Ext_c(K,D) \] 
and 
\[ \check H^2(X,\uline Z)  
\cong  \check H^2(X,\uline{Z_0}) \oplus  \check H^2(X,\uline{D}) 
\cong  \check H^3(X,\Gamma) \oplus  \check H^2(X,D).\] 
This splits the problem to determine the obstruction 
class into two cases, where $Z$ is either discrete or 
a quotient $\fz/\Gamma$. 
\end{rem} 

\begin{defn} \mlabel{def:kgrp} We call a connected abelian topological group 
$Z$ a {\it $K(\Gamma,1)$-group} if 
it is locally contractible, $\pi_1(Z) \cong \Gamma$, and $\tilde Z$ is contractible 
and divisible.\begin{footnote}{$K(\Gamma,1)$ groups exist for every 
abelian group $\Gamma$. Let $G(\Gamma)$ be the group of all measurable functions 
$[0,1] \to \Gamma$ with 
finitely many values, modulo functions supported in a zero set and  
endowed with the metric $d(f,g) := |\{f \not= g\}|$, 
then $G(\Gamma)$ is contractible, locally contractible and contains 
the subgroup $\Gamma$ of constant functions as a discrete subgroup. 
If $\Gamma$ is divisible, then $G(\Gamma)$ inherits  
this property. Since every abelian group $\Gamma$ 
can be embedded in a divisible group $D$, the group $G(D)/\Gamma$ 
is a $K(\Gamma,1)$-group.
}\end{footnote}
Typical examples arise for $Z = \fz/\Gamma$, 
where $\fz$ is a topological vector space and $\Gamma \subeq \fz$ is a discrete 
subgroup. 
\end{defn}

The following examples show how obstruction classes are connected 
to various well known constructions in topology and group theory. 

\begin{ex}  \mlabel{ex:2.6} (a) (Chern classes) Let $Z$ be a $K(\Gamma,1)$-group. 
Then $\tilde Z$ is a central extension 
of $Z$ by $\Gamma$, and the obstruction class 
defines a homomorphism 
\[ \obs^{\tilde Z} = \iota_1 \: \check H^1(X,\uline{Z}) \to 
\check H^2(X,\uline{\Gamma}) 
\cong \check H^2(X,\Gamma) \cong  H^2_{\rm sing}(X,\Gamma),\] 
which is a group isomorphism assigning to a $Z$-bundle $P \to X$ its 
{\it Chern class} $\Ch(P) \in H^2(X,\Gamma)$.

(b) If $K$ is discrete, then $BK$ can be realized as a CW complex
(\cite[Thm.\ 5.1.15]{Ros94}), hence in particular as a locally contractible space. 
Now $\check H^3(BK,\Gamma) \cong H^3(BK,\Gamma) \cong 
H^3_{{\rm grp}}(K,\Gamma)$ is the group cohomology of 
$K$ with values in the trivial $K$-module $\Gamma$ 
(\cite{EML45}). 
For the universal $K$-bundle $q \: EK \to BK$, 
we obtain for a $K(\Gamma,1)$-group $Z$ the obstruction map 
\begin{multline*}
\iota_3 \circ \obs_{EK} \: \Ext_c(K,Z) \cong H^2_{{\rm grp}}(K,Z)
 \to \check{H}^3(BK,\Gamma) \cong H^3_{{\rm grp}}(K,\Gamma). 
\end{multline*}
It is easy to see that this coincides with the  natural connecting map 
\[ \delta_2 \: H_{{\rm grp}}^2(K,Z) \to H_{{\rm grp}}^3(K,\Gamma) \] 
in the long exact sequence in 
group cohomology induced by the short exact sequence $\Gamma \into \tilde Z
\onto Z$.

(c) Let $q_X \: \tilde X \to X$ be a simply connected covering 
of $X$, 
consider the group $K = \pi_1(X)$ of deck transformations as a discrete group,  
and $\tilde X$ as a principal $K$-bundle. Then 
the obstruction class leads to a homomorphism 
\[ \obs_{\tilde X} \: H^2_{{\rm grp}}(\pi_1(X),\Gamma) 
 \to \check H^2(X,\Gamma)\cong H^2(X,\Gamma) \] 
for every discrete abelian group $\Gamma$. 
Theorem~\ref{thm:7.7} below implies that $\obs_{\tilde X}$ is injective 
and that its range is the group 
$\Lambda^2(X,\Gamma) := \ker \alpha_2$ of aspherical 
cohomology classes. This is an obstruction theoretic interpretation of 
Hopf's Theorem asserting that $\Lambda^2(X,D) \cong H^2_{\rm grp}(\pi_1(X),D)$ 
(\cite{EML45}). 
\end{ex}

\section{Topological data of central Lie group extensions} \mlabel{sec:liegrp}

To understand the homomorphism $\obs_P$ 
for connected $K$, it is important to know which topological data 
coming along with a central $Z$-extension $\hat K$ of $K$ can be used 
to express $\obs_P$ more directly. Clearly, the homomorphisms  
$\partial_2^{\hat K}$ and $\partial_1^{\hat K}$ are two such pieces of data 
appearing in Theorem~\ref{thm:main}. To see what else we have, 
it is instructive to take a look at those extensions for which $\partial_2^{\hat K}$ 
vanishes. 

First we note that, for 
connected $K$ and every abelian topological group $Z$, we have a natural 
homomorphism 
\[  E^{*}\: \Hom(\pi_{1}(K), Z)\to \Ext_{c}(K,Z),\quad \gamma\mapsto 
[\gamma_*\tilde K], \] 
where 
\[ \gamma_*\tilde K :=  (\tilde K \times Z)/\{ (d, -\gamma(d)) \: 
d \in \pi_1(K)\} \]  
is the central extension of $K$ by $Z$ associated to the universal 
covering $q_K \: \tilde K \to K$ with kernel $\pi_1(K)$ by the homomorphism 
$\gamma$. 

\begin{rem} \mlabel{rem:3.1a} 
Suppose that $K$ is connected.

(a) If $D$ is a discrete abelian group, 
then \cite[Prop.~2.6]{Ne02} implies that the homomorphism 
\[ E^* \: \Hom(\pi_1(K),D) \to \Ext_c(K,D) \] 
is an isomorphism of abelian groups. 
For the so obtained central extensions $\hat K = \gamma_*\tilde K$ 
we have $\partial_2^{\hat K} = 0$ and $\partial_1^{\hat K} = \gamma$, so that 
all information is encoded in~$\partial_1^{\hat K}$. 

(b) If $Z$ is a $K(\Gamma,1)$-group and 
$\hat K$ a central extension of $K$ by $Z$, 
then the Chern class $\Ch([\hat K])$ is an 
element of $H^2_{\rm sing}(K,\Gamma)$. 
From the Universal Coefficient Theorem we 
obtain the short exact sequence
\[ 0 \to \Ext(H_1(K), \Gamma) \cong \Ext(\pi_1(K), \Gamma) 
\to H^2(K,\Gamma) \to 
\Hom(H_2(K),\Gamma) \to 0,\]
where we have used that $\pi_1(K)$ is abelian, so that 
$H_1(K) \cong \pi_1(K)$ by the Hurewicz Theorem.  
The $Z$-bundles over $K$ corresponding to extensions of the form 
$\gamma_* \tilde K$ for $\gamma \in \Hom(\pi_1(K),Z)$ 
define the trivial homomorphism $H_2(K) \to \Gamma$ 
which can be interpreted as 
the ``curvature'' of $\hat K$. 
Since the $Z$-bundle $\gamma_*\tilde K$ is topologically trivial 
if and only if the homomorphism $\gamma \: \pi_1(K) \to Z$ 
lifts to a homomorphism 
$\tilde\gamma \: \pi_1(K) \to \tilde Z$, we see that 
$\im(E^{*}) \cong \Ext(\pi_1(K),\Gamma)$, and we conclude that the 
Chern class provides a surjection 
\[ \Ch : \Ext_c(K,Z)_{\rm flat} := \im(E^*) 
\to \Ext(\pi_1(K),\Gamma),\] 
showing that all flat $Z$-bundles over $K$ can be realized by 
flat central $Z$-extensions. 
\end{rem}

When specialized to Lie groups and $Z = \fz/\Gamma$, the 
following proposition corrects a wrong claim in \cite[Cor.~7.15]{Ne02}.

\begin{prop} \mlabel{prop:corrne02-top} 
Let $Z$ be a $K(\Gamma,1)$-group and $K$ be a connected locally contractible 
topological group. 
For the subgroup 
\[ \Ext_c(K,Z)_0 := \{ [\hat K] \in \Ext_c(K,Z) \: \partial_2^{\hat K} = 0\}, \] 
we have an exact sequence 
\[ \Ext_c(K,\tilde Z) \smapright{(q_Z)_*} 
\Ext_c(K,Z)_0 \sssmapright{\zeta} \Ext(\pi_1(K),\Gamma), \quad 
\zeta([\hat K]) = [\pi_1(\hat K)].\] 
\end{prop} 

Note that $\pi_1(\hat K)$ is abelian because $\hat K$ is a topological group. 

\begin{prf} If $K^\sharp$ is a central $\tilde Z$-extension of $K$, 
then $\hat K := (q_Z)_*K^\sharp = K^\sharp/\Gamma$ 
is a central extension for which 
$\partial_2^{\hat K} = \pi_1(q_Z) \circ \partial_2^{K^\sharp} = 0$ 
follows from $\pi_1(q_Z) = 0$. 
Moreover, 
the isomorphism $\pi_1(K^\sharp) \to \pi_1(K)$ and the natural 
homomorphism $\pi_1(K^\sharp)\to \pi_1(\hat K)$ lead to a splitting of 
the extension $\pi_1(\hat K)$ of $\pi_1(K)$ by $\Gamma$, so that 
$[\hat K] \in \ker \zeta$. 

If, conversely, $[\hat K] \in \ker \zeta$, then 
Proposition~\ref{prop:8.7}(ii) implies that $\hat K$ is a 
trivial $Z$-bundle, hence of the form 
$\hat K = Z \times K$ with a multiplication 
\[ (z,k)(z',k') = (zz'f(k,k'), kk'),\] 
where $f \: K \times K \to Z$ is a continuous $2$-cocycle. 
The splitting homomorphism $\sigma \: \pi_1(K) \to \pi_1(\hat K) \subeq 
(\hat K)\,\tilde{}$ maps onto a discrete central subgroup of 
$(\hat K)\,\tilde{}$. The latter group is of the form 
$\tilde Z \times \tilde K$ with a multiplication of the form 
\[ (z,k)(z',k') = (zz'\tilde f(k,k'), kk'),\] 
where $\tilde f \: \tilde K \times \tilde K \to \tilde Z$ 
is a continuous $2$-cocycle satisfying 
$q_z \circ \tilde f = f \circ (q_K \times q_K)$. 
We conclude that $K^\sharp := (\hat K)\,\tilde{}/\im(\sigma)$ is a 
central $\tilde Z$-extension of $K$ with 
$K \cong K^\sharp/\Gamma = (q_Z)_*K^\sharp$. 
This completes the proof. 
\end{prf}

\begin{prop}
  \mlabel{prop:0.1} 
Suppose that $\hat K$ is an extension of the connected 
group $K$ by the $K(\Gamma,1)$-group $Z$. 
If 
\[ \partial_2^{\hat K} \: \pi_2(K) \to \Gamma \quad \mbox { and } \quad 
[\pi_1(\hat K)] \in \Ext(\pi_1(K),\Gamma) \] 
vanish, then $\obs^{\hat K}([P])$ vanishes for every $K$-bundle~$P$. 
\end{prop}

\begin{prf} Proposition~\ref{prop:corrne02-top} implies the existence 
of a central extension 
$K^\sharp$ of $K$ by $\tilde Z$ with $K^\sharp/\Gamma \cong \hat K$. 
Therefore the obstruction class 
$\obs^{\hat K}([P]) \in \check H^2(X,\uline{Z})$ is the 
image of the corresponding obstruction class 
$\obs^{K^\sharp}([P]) \in \check H^2(X,\uline{\tilde Z})=\{0\}$ 
(Lemma~\ref{lem:0.2}(c)), hence trivial. 
\end{prf}

\begin{rem} \mlabel{rem:3.4b} 
Since $\tilde Z$-extensions of $K$ define trivial obstruction classes,
the preceding proof shows that the information in $[\hat K]$ that 
is relevant for $\obs_P([\hat K])$ only depends on 
$[\hat K]$ modulo the image of $(q_Z)_*$. 
Two such classes coincide if and only if 
\[ \partial_2^{\hat K_1} = \partial_2^{\hat K_2} \: \pi_2(K) \to \Gamma 
\quad \mbox{ and } \quad 
0 = \zeta([\hat K_1] - [\hat K_2]) \in \Ext(\pi_1(K),\Gamma).\] 
\end{rem}

\begin{rem} \mlabel{rem:0.4} 
In view of Proposition~\ref{prop:corrne02-top} and 
Remark~\ref{rem:3.1a}(b),  
\begin{equation} \label{eq:decompo}
\Ext_c(K,Z)_0 = \im((q_Z)_*) + \im(E^*).
\end{equation}

Flat central extensions are precisely those whose Chern class 
defines the trivial homomorphisms $H_2(K) \to \Gamma$ 
and $\Ext_c(K,Z)_0$ consists of those central extensions for which 
this homomorphism vanishes on the image of $\pi_2(K)$ under the 
Hurewicz homomorphism $h_2 \: \pi_2(K) \to H_2(K)$. 

Let $[\hat K] \in \Ext_c(K,Z)_0$. 
In view of \eqref{eq:decompo}, $[\hat K]$ is a sum of a flat 
extension and an extension of the form $(q_Z)_*K^\sharp$, where  
$K^\sharp$ is a central extension of $K$ by $\tilde Z$. 
Since the $Z$-bundle $(q_Z)_*K^\sharp$ 
over $K$ has a trivial Chern class, 
it follows that, whenever $\partial_2^{\hat K}$ vanishes, then also 
the homomorphism $H_2(K) \to \Gamma$ defined by the Chern class of 
the $Z$-bundle $\hat K$ vanishes. Therefore no additional information 
is gained by considering $H_2(K)$ instead of $\pi_2(K)$. 
\end{rem}

\section{\v{C}ech cohomology on spheres} 
\mlabel{sec:sn}

We introduce a good cover 
of $\bS^n$ with respect to which we specify \v{C}ech cohomology classes, 
the isomorphisms we need in the following, 
and the suspension homomorphism in  \v{C}ech cohomology. 
We also recall the gluing construction of $K$-bundles 
over $\bS^n$ from elements of $\pi_{n-1}(K)$. 

\subsection{Calculating \v Cech cohomology classes on spheres} 
\mlabel{sect:covering}

Suppose that $Z$ is connected. 
For the spheres $\bS^n$, $n > 1$, we obtain from 
\eqref{eq:iota} 
isomorphisms  
\[ \iota_{n-1} \: \check H^{n-1}(\bS^n,\uline Z)\to \check H^{n}(\bS^n,\Gamma) 
\cong \Gamma.\]
Below we have to deal with the suspension map 
$\susp{1}\:\check H^{1}(\bS^2,\uline Z) 
\to \check H^{2}(\bS^3,\uline Z)$. To make this map more explicit, we 
shall identify both sides with $\Gamma$, and to do that, we need  
explicit isomorphisms $S_n \: \check H^n(\bS^n,\Gamma) \to \Gamma$. 

To this end, we realize $\bS^n$ as the (relative) boundary of the standard simplex 
\[ \Delta^{n+1} = [e_0, \ldots, e_{n+1}] \subeq \R^{n+2},\]  
where $e_0,\ldots, e_{n+1}$ are the canonical basis vectors in 
$\R^{n+2}$. Then the $j$th face 
\[ F_j := [e_0, \ldots, \hat{e_j}, \ldots, e_{n+1}], \qquad 
j=0,\ldots, n+1 \] 
is a closed subset of $\bS^n = \partial \Delta^{n+1}$ and we 
write 
$\Lambda_j := F_j^c$ for its open complement, which is the intersection 
of $\partial\Delta^{n+1}$ with a convex subset of $\Delta^{n+1}$. The 
latter property is inherited by all finite intersections of the 
$\Lambda_j$, which immediately implies that 
$\cU := (\Lambda_0,\ldots, \Lambda_{n+1})$ is a {\it good cover} of $\bS^n$, 
i.e., all intersections of covering sets are contractible. 
Hence we can calculate \v Cech cohomology directly with the cover 
$\cU$ 
via\begin{footnote}{From the isomorphisms 
$\check H^n(X,\uline Z)\cong \check H^{n+1}(X,\Gamma)$, 
it follows that these cohomology groups vanish if $X$ 
is contractible. Hence the claim	 follows from \cite[Cor.,p.~213]{God}.}
\end{footnote}
\[ \check H^k(\bS^n, \uline Z) \cong \check H^k_{\cU}(\bS^n, \uline Z)
\quad \mbox{ and } \quad 
 \check H^k(\bS^n, \Gamma) \cong \check H^k_{\cU}(\bS^n, \Gamma).\]    
We also note that $\bigcup_{j=0}^{n+1} F_j = \bS^n$ implies that 
$\bigcap_{j = 0}^{n+1} \Lambda_j = \eset$. 
In particular, the space $\check C^{n+1}_{\cU}(\bS^n,\Gamma)$ 
of $(n+1)$-cochains is trivial, which implies that 
\[ \check H^n_{\cU}(\bS^n,\Gamma) 
\cong C^n_{\cU}(\bS^n,\Gamma)/\delta(C^{n-1}_{\cU}(\bS^n,\Gamma)).\] 
Elements of $\check C^n_{\cU}(\bS^n,\Gamma)$ assign to 
$(n+1)$-fold intersections 
\[ \Lambda_{i_1,\ldots, i_{n+1}} 
:= \bigcap_{j = 1}^{n+1} \Lambda_{i_{j}} \] 
elements of $\Gamma$. Since there are only $n+2$ such intersections: 
\[ \Lambda_{(j)} := \Lambda_{0,1,\ldots, \hat{j}, \ldots, n+1} 
= \bigcap_{k\not=j} \Lambda_k=\inn(F_{j})\] 
we have $C^n_{\cU}(\bS^n,\Gamma) \cong \Gamma^{n+2}$. 

\begin{lem} \mlabel{lem:iso-sn} The map 
\begin{equation}
  \label{eq:sn}
\tilde S_n \:  C^{n}_{\cU}(\bS^n, \Gamma) \cong \Gamma^{n+2} \to \Gamma, 
\quad \tilde S_n(\gamma_0,\ldots, \gamma_{n+1}) 
:= \sum_{j = 0}^{n+1} (-1)^j \gamma_j, 
\end{equation}
factors through an isomorphism 
$S_n \: \check H^n(\bS^n,\Gamma) \to \Gamma.$
\end{lem}

\begin{prf} 
For $\beta \in C^{n-1}_{\cU}(\bS^n,\Gamma)$, the coboundary is given by 
\begin{align*}
\delta(\beta)_{(j)} 
&= \beta_{1,\ldots, \hat j, \ldots, n+1} 
- \beta_{0,2,\ldots, \hat j, \ldots, n+1} \pm \cdots \\
&+(-1)^{j-1} \beta_{0,1,\ldots,j-2,\hat{j-1}, \hat j, \ldots, n+1} 
+(-1)^j \beta_{0,1,\ldots, \hat j, \hat{j+1},\ldots, n+1} 
\pm \cdots.
\end{align*}  
For the cochains $\beta^k$, $k = 0,\ldots, n$, assigning the value $\gamma$ to 
$\Lambda_{0,1,\ldots, \hat{k},\hat{k+1},\ldots, n+1}$ 
and $0$ to all other $(n-1)$-fold intersections, we obtain
\[ \delta(\beta^k)_{(j)} 
= (-1)^{j-1} \delta_{j-1,k} \gamma 
+ (-1)^j \delta_{j,k} \gamma 
= (-1)^k(\delta_{j-1,k} + \delta_{j,k}) \gamma.\] 
This implies
$\delta(C^{n-1}_{\cU}(\bS^n, \Gamma))\subseteq \ker(\tilde S_n)$. On the other hand, if
$(\gamma_{0},\ldots,\gamma_{n+1})$ is contained in $\ker(\tilde S_n)$, then $\gamma_{n+1}=(-1)^{n}\sum_{j=0}^{n}(-1)^{j}\gamma_{j}$ implies
\begin{align*}
 (\gamma_{0},\ldots,\gamma_{n+1})=(\gamma_{0},\gamma_{0},0,\ldots,0)&-
(0,\gamma_{0}-\gamma_{1},\gamma_{0}-\gamma_{1},0,\ldots,0)\pm\ldots\\ 
\ldots&+(-1)^{n}\Big(0,\ldots,0,\sum_{j=0}^{n}(-1)^{j}\gamma_{j},\sum_{j=0}^{n}(-1)^{j}\gamma_{j}\Big).
\end{align*}
Since the right-hand side is contained in $\delta(C^{n-1}_{\cU}(\bS^n, \Gamma))$, this shows that indeed $\delta(C^{n-1}_{\cU}(\bS^n, \Gamma))=\ker(\tilde S_n)$.
\end{prf}

\begin{defn}\mlabel{defalpha}
Since $\Gamma$ is discrete, Lemma~\ref{lem:iso-sn} provides 
for each $n > 0$ a natural isomorphism 
$S_n \: \check H^n(\bS^n, \Gamma) \to \Gamma$. 
For a paracompact locally contractible space $X$, 
we thus obtain natural homomorphisms 
\begin{equation}
  \label{eq:alphan}
\alpha_n \: 
\check H^n(X,\Gamma) \to \Hom(\pi_n(X),\Gamma),
\quad \alpha_n(c)([\sigma]) = S_n(\sigma^*c)
\end{equation}
whose kernel $\Lambda^{3}(X,\Gamma)$ we call the \emph{aspherical classes}.

In fact, for  the natural isomorphism 
$\check H^n(X,\Gamma) \to H^n(X,\Gamma)$ 
(\cite[p.~184]{Br97b}), the evaluation map 
\[\beta_n \: H^n(X,\Gamma) \to \Hom(H_n(X),\Gamma) \] 
coming from the Universal Coefficient Theorem 
(cf.\ \cite{Br97a}) 
and the Hurewicz homomorphism 
\[ h_n \: \pi_n(X) \to H_n(X), \] 
we have
\begin{equation}\label{eq:hurw-rel}
\alpha_n(c) = \beta_n(c) \circ h_n.
\end{equation}
Since $\alpha_n$, $\beta_n$ and $h_n$ are functorial in $X$, 
it suffices to verify this relation for $X = \bS^n$, 
where it is an immediate consequence of the compatibility 
of both sides with the suspension map 
(cf.\ Proposition~\ref{Dreiecksdiagramm} below and 
the Mayer--Vietoris Sequence in singular (co)homology 
\cite[\S IV.15]{Br97a}). 
\end{defn}

\subsection{The suspension map in \v Cech  cohomology} 

For each abelian topological group $A$, there is a suspension homomorphism 
\begin{equation*}
\susp{k}:\check H^k(\bS^n,\underline{A}) \to \check H^{k+1}(\bS^{n+1},
\underline{A}). 
\end{equation*}
To describe this map, observe that the good open cover we used 
for $\bS^{n+1}$ reduces to the one defined for $\bS^n$ by restriction to the 
equator $\partial F_{n+2}\cong \bS^{n}\subseteq \R^{n+1}$. 
Moreover, there exists a retraction 
\begin{equation}\label{eqn:susp1}
p\: \bS^{n+1}\backslash (\inn(F_{n+2})\cup\{e_{n+2}\})\to \partial F_{n+2} \cong \bS^n,
\end{equation}
restricting for $0 \leq j \leq n+1$ to a 
retraction $F_{j}\backslash \{e_{n+2}\} \to F_{j}\cap F_{n+2}$. 
It satisfies $p(\Lambda_{j}\backslash \inn(F_{n+2}))
\subseteq \Lambda_{j}\cap \R^{n+1}$ and thus maps the open subset 
$\Lambda_{j}$ of $\bS^{n+1}$ to the corresponding 
open subset $\Lambda_{j} \cap \bS^{n}$. 

The suspension map $\susp{k}$
is defined by sending
a cochain $(g_{i_0,\ldots,i_{k}})$ on $\bS^n$ to the cochain 
$(h_{i_0,\ldots,i_{k+1}})$ on $\bS^{n+1}$ defined by
\[h_{i_0,\ldots,i_{k+1}}:=
\begin{cases}
g_{i_0,\ldots,i_{k}}\circ p & \text{ for } i_0 < \ldots < i_{k+1} = n+2 \\ 
0 & \text{ for } n+2 \not\in \{ i_0,\ldots,i_{k+1}\}. 
\end{cases}\]
We work in this article always with totally 
antisymmetric \v{C}ech cocycles (\cite[Rem.~4.01, p.~218]{Ha77}), 
so if $i_{n+1}$ appears in the list of indices, it may always be moved 
to the last position. The so 
defined cochain $(h_{i_0,\ldots,i_{k+1}})$ 
then satisfies 
\[ (\delta h)_{i_0,\ldots,i_{k+2}}=
\begin{cases}
(\delta g \circ p)_{i_0,\ldots,i_{k+1}}& \text{ for } i_{k+2} = n+2 \\ 
0 & \text{ otherwise}.   
\end{cases}\]
This shows that, if $g$ is a cocycle, then so is $h$, and we likewise 
see that, if $g$ is a coboundary, then $h$ is a coboundary.

\begin{prop}
\mlabel{Dreiecksdiagramm}
If $A=\Gamma$ is discrete, then the diagram
\[ 
\xymatrix@=1.5em{\check H^n(\bS^n,\Gamma)\ar[dr]^(.6){S_{n}} \ar[rr]^(.45){\susp{n}}&&\check H^{n+1}(\bS^{n+1},\Gamma)\ar[dl]_(.625){S_{n+1}}\\&\Gamma}
\]
commutes and $\susp{n}$ is an isomorphism. If $\Gamma$ is discrete in the topological vector space
$\mathfrak{z}$ and  $Z=\mathfrak{z}/\Gamma$, then the diagram
\begin{equation*}
 \xymatrix@=1.5em{
\check{H}^{k}(\bS^{n},\uline{Z})\ar[rr]^(.45){\susp{k}}\ar[d]^{\iota_{k}}&&\check{H}^{k+1}(\bS^{n+1},\uline{Z})\ar[d]^{\iota_{k+1}}\\
\check{H}^{k+1}(\bS^{n},\Gamma)\ar[rr]^(.45){\susp{k+1}}&&\check{H}^{k+2}(\bS^{n+1},\Gamma)}
\end{equation*}
also commutes.
\end{prop}

\begin{prf}
 From the construction of $S_{n},S_{n+1}$ and $\susp{n}$ one immediately deduces 
\break $S_{n+1} \circ \susp{n} = S_n$.
 This also implies that $\susp{n}$ is an isomorphism since $S_{n}$ and $S_{n+1}$ are so.

 The isomorphism $\iota_{k}$ is defined by sending the cocycle $(g_{i_{0},\ldots,i_{k}})$ to the cocycle $\delta (\tilde{g})_{i_{0},\ldots,i_{k+1}}$, 
 where $\tilde{g}_{i_{0},\ldots,i_{k}}$ are arbitrarily chosen lifts of $g_{i_{0},\ldots,i_{k}}$.
 Since \break $(\tilde{g}_{i_{0},\ldots,i_{k}}\circ p)$ lifts the image $(g_{i_{0},\ldots,i_{k}}\circ p)$
 of $g_{i_{0},\ldots,i_{k}}$ under $\susp{k}$, the commutativity of the second diagram 
follows.
 \end{prf}

The above description also carries over to the
case of a non-abelian topological group $K$. In this case it defines a map 
\begin{equation*}
\susp{0} \: \check{H}^{0}(\bS^{n},\underline{K})\to \check{H}^{1}(\bS^{n+1},\underline{K}), 
\quad 
\susp{0}(g)_{ij} = 
\begin{cases}
g_i \circ p  & \text{for } j = n + 2 \\
\1  & \text{for } j <  n + 2.   
\end{cases}
\end{equation*}

\subsection{Bundles over spheres} 
\mlabel{subsec:2.3}

\begin{defn}
For connected $K$, elements of $\Bun(\bS^n,K)$ are classified by 
homotopy classes $[f] \in [\bS^n, BK] \cong \pi_n(BK)$. Since
\[ \partial^{EK}_n \: \pi_n(BK) \to \pi_{n-1}(K) \] 
is an isomorphism for $n \geq 1$, the characteristic class 
of $P \cong f^*EK$ is given by 
\begin{equation}
  \label{eq:defh}
h([P]) := \partial^{EK}_n([f]) = \partial^P_n([\id_{\bS^n}]) 
\in \pi_{n-1}(K),
\end{equation}
where the second equality follows from Lemma~\ref{lem:0}.
\end{defn}

\begin{rem} \mlabel{rem:3.4} (a) 
Applying Lemma~\ref{lem:0} to $Y = \bS^n$, we obtain 
for a $K$-bundle $P \to X$ 
the following interpretation of $\partial_n^P$ in terms 
of characteristic classes of bundles over $\bS^n$: 
\begin{equation}
  \label{eq:charac}
\partial^P_n([f]) = \partial^{f^*P}_n([\id_{\bS^n}]) 
=  h(f^*P) \in  \pi_{n-1}(K).
\end{equation}

(b) For $n = 2$ and $Z = \fz/\Gamma$, 
we obtain in particular for any $Z$-bundle $Q$ over~$\bS^2$ 
\begin{equation}
  \label{eq:charac2}
h(Q) = \partial^Q_2([\id_{\bS^2}]) \in  \pi_{1}(Z)\cong \Gamma.
\end{equation}
\end{rem}

\begin{rem} \mlabel{rem:3.1}
Let $n \geq 2$ and
$f \: \bS^{n-1} \to K$ be a continuous map. Below we recall the clutching 
construction of a principal $K$-bundle $P_f \to \bS^n$ satisfying
$\partial_n^{P_f}([\id]) = [f]$ and $[P_{f}]=\susp{0}(f)$, 
where $\partial_n^{P_f}$ is the connecting map 
from \eqref{eq:conmap} and $\susp{0}$ is the suspension map.

To this end, 
let $\bB^n$ be the closed unit ball in $\R^n$ and 
$\eta_j \: \bB^n \to \bS^n$, $j = 0,1$, be topological embeddings with 
$\eta_0(\bS^{n-1}) = \eta_1(\bS^{n-1})$ and for which the 
closed subsets $A_j := \im(\eta_j)$ satisfy 
\[ \bS^n = A_0 \cup A_1 \quad \mbox{ and } \quad 
A_0 \cap A_1 = \eta_j(\bS^{n-1}) \cong \bS^{n-1}.\] 

We now consider the quotient 
\[ P_f := (\bB^n \times \{0,1\} \times K)/\sim,\] 
where the equivalence relation is given by 
$(x,j,k) \sim (x',j',k')$ if and only if 
$x = x'$, $j = j'$, and $k = k'$, or 
$x \in \bS^{n-1}$, $\eta_0(x) = \eta_1(x)$ 
and $k = f(x)k'$. We write $[(x,j,k)]$ for the equivalence class
of $(x,j,k)$ in $P_f$. Then it 
is easy to see that $K$ acts from the right on $P_f$ by 
$[(x,j,k)].k' := [(x,j,kk')]$ and 
that $q \: P_f \to \bS^n, [(x,j,k)] \mapsto \eta_j(x)$ 
is a continuous map whose fibers are the $K$-orbits. 
Clearly, 
$\sigma_j(\eta_j(x)) := [(x,j,\1)]$, $j =1,2$, define continuous sections on the 
closed subsets $A_j$, and on $A_0 \cap A_1$  
these sections satisfy 
\[ \sigma_1(\eta_1(x)) = \sigma_0(\eta_0(x)).f(x) 
\quad \mbox{ for } \quad x \in \partial \bS^{n-1}.\] 
Extending the map $\eta_0(x) \mapsto f(x)$ to a neighborhood of 
$A_0 \cap A_1$ in $A_0$, the assignment 
$x \mapsto \sigma_0(\eta_0(x)).f(x)$ 
defines a continuous section on an open 
neighborhood of $A_1$ extending $\sigma_1$. This implies 
that $P_f$ is locally trivial. 

To verify the relation 
\begin{equation}
  \label{eq:charclsphere}
\partial_n^{P_f}([\id_{\bS^n}]) = [f] \in \pi_{n-1}(K), 
\end{equation}
we consider the quotient map 
\[ p \: \hat A_0 := \bS^{n-1} \times I \to A_0, \quad 
(x,t) \mapsto \eta_0(tx) \] 
collapsing $\bS^{n-1} \times \{0\}$ to a single point. 
Gluing $\hat A_0$ in the canonical way to $A_1$, we obtain a 
space $\hat \bS^n := \hat A_0 \cup A_1$ homeomorphic to $\bB^n$ and a map 
$\xi \: \hat\bS^n \to \bS^n$ inducing a homeomorphism 
$\oline\xi \: \hat \bS^n/(\bS^{n-1} \times \{0\}) \to \bS^n$. 
Now 
\[ \tilde\xi \: \hat\bS^n \to P_f, \quad 
\tilde\xi(x,t) \mapsto \sigma_0(\eta_0(tx)).f(x) \] 
defines a continuous lift of $\oline\xi$. 
Its restriction to $\bS^{n-1} \times \{0\}$ 
is given by $(x,0) \mapsto \sigma_0(\eta_0(0)).f(x)$, and this implies 
that $\partial_n^{P_f}([\id]) = [f]$.

In order to verify the relation $[P_{f}]=\susp{0}(f)$, 
we observe that we may identify $\bS^{n}$ with
$\partial \Delta^{n+1}$ in such a way that $F_{n+2}$ gets identified with
$A_{2}$ and the closure of ${\Lambda_{n+1}}$ with $A_{1}$. Then we can use the 
map $f \circ p$ to extend the section $\sigma_{2}$ to $A_{1}\backslash \{e_{n+2}\}$
and the transition functions of this system of sections (restricted to the $\Lambda_{i}$'s)
has precisely $\susp{0}(f)$ as associated \v{C}ech cocycle.
\end{rem}

We end this section with showing that suspension and clutching are
compatible in the following sense. This will be the basis for our
upcoming arguments.

\begin{prop}\mlabel{Quadratdiagramm}
 For an arbitrary abelian topological group $Z$, 
$[\hat{K}]\in \Ext_{c}(K,Z)$ and $n\in\N$
 the diagram
\begin{equation}\label{eqn:Quadratdiagramm}
	\vcenter{
 	\xymatrix{\
	\check{H}^0(\bS^n,\underline{K})  \cong C(\bS^n,K)
		\ar[d]^{\delta_0^{\bS^n}}
		\ar[rr]^{\susp{0}}_{f\mapsto[P_{f}]}
	&&
	\check{H}^1(\bS^{n+1},\underline{K}) 
		\ar[d]^{\delta_1^{\bS^{n+1}}}
	\\
	\check{H}^1(\bS^n,\underline{Z}) 
	\ar[rr]^{-\susp{1}}
	&&
	\check{H}^2(\bS^{n+1},\underline{Z})
	}}
\end{equation}
commutes. Moreover, we have $\delta_{0}^{\bS^{n}}(f)=[f^{*}\hat{K}]$.
\end{prop}

\begin{prf} We first perform some basic constructions. 
Recall the covering $\Lambda_{0},\ldots, \Lambda_{n+2}$
of $\bS^{n+1}\cong \partial \Delta^{n+2}$ from Section \ref{sect:covering}.
\begin{enumerate}
	\item For $\tilde{f}\: \bS^{n}\cong \partial F_{n+2}\to K$ we set
	\begin{equation*}
	 f:=\tilde{f} \circ p\: \bS^{n+1}\backslash (\inn(F_{n+2})\cup \{e_{n+2}\})\to K,
	\end{equation*}
	where $p$ is the projection map from \eqref{eqn:susp1}. With this we define
	\begin{equation*}
	 f_{i}:=\left.f\right|_{\Lambda_{i}\backslash\inn(F_{n+2})}\to K. 
	\end{equation*}
	\item Define $g_{ij}\: \Lambda_{ij}\to K$ for $j=n+2$ by 
$g_{i\,n+2}:=f_{i}$ and $g_{ij}=\1$ if
	$n+2\notin\{i,j\}$. By the previous remark, these are transition functions 
for the bundle $P_{f}$ on the cover $(\Lambda_{i})_{i=0, \ldots, n+2}$.
	\item Since $\Lambda_{i}\cap \partial F_{n+2}$ 
is contractible for $i<n+2$, 
 there exist  lifts\break 
	$\tilde{f}_{i}\: \Lambda_{i}\cap \partial F_{n+2}\to \hat{K}$. Then $\hat{f}_{i}:=\tilde{f}_{i}\circ p$ is a lift of $f_{i}$ satisfying
	$\hat{f}_{i}=\left.\tilde{f}_{i}\right|_{\Lambda_{i}\cap \partial F_{n+2}}\circ p$.
	\item The maps $\hat{g}_{i\,n+2}:=\hat{f}_{i}$ if $i<n+2$ and $\hat{g}_{ij}=1$ if $n+2\notin\{i,j\}$ define a lifting cochain of the cocycle $g_{ij}$.
	Thus we have
	\begin{alignat*}{3}
	 \delta_{1}^{\bS^{n+1}}(g)_{ij\,n+2}&=\hat{g}_{ij}\cdot \hat{g}_{j\,n+2}
\cdot \hat{g}_{i, n+2}^{-1}=\hat{f}_{j}\cdot\hat{f}_{i}^{-1}&\quad&\text{ if }i,j<n+2\\
	 \delta_{1}^{\bS^{n+1}}(g)_{ijk}&=0&\quad&\text{ if }n+2\notin\{i,j,k\}.
	\end{alignat*}
\end{enumerate}
Since the cocycle $(g_{ij})$ represents $[P_{f}]$ and $\hat{g}_{ij}$ lifts $g_{ij}$, 
the above cocycle represents $\delta_{1}^{\bS^{n+1}}([P_{f}])$. 
On the other hand,
$\hat{f}_{i}\cdot \hat{f}_{j}^{-1}\res_{\Lambda_{ij}\cap \partial F_{n+2}}$ 
is a cocycle representing $\delta_{0}^{\bS^{n}}(f)$. Thus we have
\begin{align*}
& \big(\susp{1}(\delta_{0}^{\bS^{n}}(f))\big)_{ij\,n+2}
=(\delta_{0}^{\bS^{n}}(f))_{ij}\circ p \\
&=\Big(\hat{f}_{i}\res_{\Lambda_{ij}\cap \partial F_{n+2}}\circ p\Big) \cdot
 \Big(\hat{f}_{j}^{-1}\res_{\Lambda_{ij}\cap \partial F_{n+2}}\circ p\Big) 
=\hat{f}_{i}\cdot \hat{f}_{j}^{-1},
\end{align*}
showing our first claim. Since $\hat{f}_{i}$ 
can be used to construct sections of $f^{*}\hat{K}$ on the cover $(\Lambda_{i})_{i=0,\ldots,n+2}$ it also follows that
$\delta_{0}^{\bS^{n}}(f)=[f^{*}\hat{K}]$.
\end{prf}

\begin{rem}\mlabel{gerbes}
 The commutativity of the diagram in the previous proposition can also be understood in the language of (smooth) bundle gerbes (cf.\ \cite{Mur}, \cite{Gomi}, \cite{SW10}).
 First we note that the diagram \eqref{eqn:Quadratdiagramm} may be rephrased in the smooth setting by considering sheaves of
 germs of smooth $K$- and $Z$-valued functions. The proof of the commutativity of this diagram then carries over literally to the  
 smooth setting.

 Bundle gerbes are defined analogously to the local definition of principal bundles, except that the cohomological
 dimension is raised by one. The role of an open trivializing cover of a smooth manifold $M$ will be taken by an arbitrary surjective 
 submersion\footnote{In the case of a (trivializing) cover $(U_{i})_{i\in I}$ of $M$, the 
 surjective submersion would be \break $\bigsqcup U_{i}\to M$.} $Y\to M$, giving rise to the $k$-fold fiber products $Y^{[k]}:=Y\times_{M}\cdots \times_{M}Y$ (replacing 
 the $k$-fold intersections of the open cover). A principal 
$Z$-bundle over $M$ is then a $Z$-valued function $z\: Y^{[2]}\to Z$ 
satisfying the cocycle condition
 \begin{equation*}
  \pi_{12}^{*}z\cdot \pi_{23}^{*}z=\pi_{13}^{*}z
 \end{equation*}
 on $Y^{[3]}$. A bundle gerbe thus consists of a principal $Z$-bundle $Q\to Y^{[2]}$ and
 an isomorphism
 \begin{equation*}
  \mu\: \pi_{12}^{*}Q\otimes \pi_{23}^{*}Q\to \pi_{13}^{*}Q
 \end{equation*} 
 of principal $Z$-bundles over $Y^{[3]}$ such that the diagram
 \begin{equation}\label{eqn:associativity}
  \vcenter{\xymatrix{\pi_{12}^{*}Q\otimes \pi_{23}^{*}Q\otimes \pi_{34}^{*}Q\ar[rr]^(.55){{\pi_{123}^{*}\mu \otimes \id}}\ar[d]^{\id\otimes\pi_{234}^{*}\mu}&&\pi_{13}^{*}Q\otimes \pi_{34}^{*}Q\ar[d]^{\pi_{134}^{*}\mu}\\\pi_{12}^{*}Q\otimes \pi_{24}^{*}Q\ar[rr]^{\pi_{124}^{*}\mu}&&\pi_{14}^{*}Q}}
 \end{equation}
 of isomorphisms of principal $Z$-bundles over $Y^{[4]}$ commutes. 

 Bundle gerbes have characteristic classes in $\check{H}^{2}(M,\uline{Z})$, which can be constructed as follows.
 One chooses an open covering $(U_{i})_{i\in I}$ such that there exist sections $\sigma_{i}\: U_{i}\to Y$ of the submersions $Y\to M$, yielding
 $\sigma_{ij}\:U_{i}\cap U_{j}\to Y^{[2]}$. Let $s_{ij}$ be a section of $Q_{ij}:=\sigma_{ij}^{*}Q$. On $U_{i}\cap U_{j}\cap U_{k}$ we thus obtain
 a section $\mu_{*}(s_{ij},s_{jk})$ of $Q_{ik}$ and thus a unique smooth function $\gamma_{ijk}\: U_{i}\cap U_{j}\cap U_{k}\to Z$ satisfying
 $\mu_{*}(s_{ij},s_{jk})=s_{ik}\cdot \gamma_{ijk}$. Using 
\eqref{eqn:associativity} one easily checks that 
$(\gamma_{ijk})$ indeed comprises
 a \v{C}ech $2$-cocycle whose cohomology class 
is independent of all choices. 

 Starting with a smooth function $f\: \bS^{n}\to Z$, we can now construct two bundle gerbes of $\bS^{n+1}$. The first one is the {\it lifting bundle gerbe},
 obtained by taking as surjective submersion the bundle $P:=P_{f}\to \bS^{n+1}$ from Remark \ref{rem:3.1} 
(which actually can be realized as a smooth principal
$K$-bundle).
Then $P^{[2]}\cong P\times K$ and we take $Q:=P\times\hat{K}\to P\times K$ as principal $Z$-bundle.
With respect to the identifications $P^{[3]}\cong P\times K\times K$ we then have
\begin{align*}
 &\pi_{12}^{*}(Q)=P\times \hat{K}\times K ,\\
 &\pi_{23}^{*}(Q)=P\times {K}\times \hat{K} \text{ and}\\
 &\pi_{13}^{*}(Q)=\{(p,k,k',\hat{k})\in P\times K\times K\times\hat{K}:kk'=q(\hat{k})\}
\end{align*}
 with the corresponding natural maps to $P\times K\times K$ (projection to the second and third entry in the last case).
 Then $\pi_{12}^{*}(Q)\otimes \pi_{12}^{*}(Q)\cong P\times (\hat{K}\times \hat{K})/Z$ (for the anti-diagonal action of $Z$) and $\mu$ is given by
 \begin{equation*}
 \mu(p,\hat{k},\hat{k'})= (p,q(\hat{k}),q(\hat{k'}),\hat{k}\hat{k'})
 \end{equation*}
The diagram \eqref{eqn:associativity} commutes for this choice of $\mu$ since the multiplication on $\hat{K}$ is associative.

 The second bundle gerbe is the {\it suspension gerbe} 
constructed by taking as surjective submersion the one obtained from an 
open cover of $\bS^{n+1}$ by two open hemispheres 
 $U_{1},U_{2}$ intersecting in
 $\left]-\varepsilon,\varepsilon\right[\times \bS^{n}$ 
(where $\bS^{n}$ is identified with the equator in $\bS^{n+1}$). 
Then $f^{*}\hat{K}$ is a principal
 $Z$-bundle over $U_{1}\cap U_{2}$ and we obtain a $Z$-bundle over
 \begin{equation*}
 Y^{[2]}:=U_{1}\sqcup (U_{1}\cap U_{2})\sqcup(U_{2}\cap U_{1})\sqcup U_{2}
 \end{equation*} 
 by taking $-f^{*}\hat{K}$ on $U_{2}\cap U_{1}$ and the trivial $Z$-bundle over $U_{1}$ and $U_{2}$.
 Since $Y^{[3]}=Y^{[2]}$ and since tensoring an arbitrary bundle with the trivial bundle yields the same bundle (up to a 
 canonical isomorphism)
 we have that $\pi_{12}^{*}P\otimes \pi_{23}^{*}P$ and $\pi_{13}^{*}P$ are canonically isomorphic. If we take this isomorphism as $\mu$ it clearly makes 
 \eqref{eqn:associativity} commute and thus completes the specification of the second bundle gerbe.

From its construction described above, 
it follows that the characteristic class of the lifting 
bundle gerbe is given by
 $\delta_{1}^{\bS^{n+1}}([P_{f}])$ and the characteristic class of the suspension bundle gerbe is given by $\susp{1}([f^{*}\hat{K}])$. Thus the
 commutativity of \eqref{eqn:Quadratdiagramm} implies that these two classes coincide (up to sign).
\end{rem}

\section{Central extensions by discrete groups} 
\mlabel{sec:1}

In Remark~\ref{rem:0.2}, we have broken up the problem 
to describe obstruction classes into the two cases where 
$Z$ is connected or discrete. In this section we take a closer look 
at the discrete case $Z = D$, i.e. where $\hat K$ is a covering of~$K$. 
Then we have to analyze the homomorphism 
\[ \obs_P \: \Ext_c(K,D) \to \check H^2(X,D).\] 
We first take a closer look at bundles 
over $\bS^2$. 

\begin{prop}\mlabel{prop:lower-degree-analog}
Let $K$ be a connected locally contractible topological group, 
$P$ be a $K$-bundle over $\bS^2$ and 
$\hat K$ be a central extension of $K$ by a discrete group 
$D$. 
Then 
\[ S_2(\delta_1([P])) = -\partial_1^{\hat K}(h(P)) \in D.\] 
\end{prop}

\begin{prf}
Let $f \: \bS^1 \cong  \partial\Delta^2 \to K$ be a 
classifying map (cf.\ Remark~\ref{rem:3.1}), so that $h([P]) = [f]$.  This
means that $P \sim P_f$, and Proposition~\ref{Quadratdiagramm} 
implies that 
\[ \susp{1}([f^*\hat K]) = - \delta_1([P_f]), \] 
so that Proposition~\ref{Dreiecksdiagramm} 
leads to 
\begin{equation*}
 S_{2}(\delta_{1}([P_{f}]))=-S_{1}([f^{*}\hat{K}]).
\end{equation*}
So it suffices to verify $S_{1}([f^{*}\hat{K}])=\partial_{1}^{\hat{K}}([f])$.
A cocycle representative of $[f^{*}\hat{K}]=\delta_{0}(f)$ 
is obtained from lifting
$\left.f\right|_{\Lambda_{i}}$ 
to $\hat{f}_{i}\:\Lambda_{i}\to\hat{K}$ for $i=0,1,2$, 
where we use the open cover $(\Lambda_{i})_{i=0,\ldots, 3}$ 
of $\partial \Delta^{2}$ from Section~\ref{sec:sn}. 
We may assume that $\hat{f}_{1}$ agrees with $\hat{f}_{2}$ on 
the open line segment $\Lambda_{12} = ]e_1, e_2[$ 
and that $\hat{f}_{2}$ agrees with $\hat{f}_{0}$ on 
$\Lambda_{02} = ]e_0, e_2[$.
Then
\begin{equation*}
 S_{1}([f^{*}\hat{K}])=\hat{f}_{1}\cdot\hat{f}_{2}^{-1}-
 \hat{f}_{0}\cdot\hat{f}_{2}^{-1}+f_{0}\cdot\hat{f}_{1}^{-1}=
 \hat{f}_{0}\cdot\hat{f}_{1}^{-1}
\end{equation*}
Note that $\hat{f}_{i}\cdot\hat{f}_{j}^{-1}$ is in $D$ and thus constant. 
Fixing a base point in $]e_0,e_1[ = \Lambda_{01}$, the 
lifts $\hat f_1,\hat f_2$ and $\hat f_0$ on 
\[ \Lambda_1 = ]e_0, e_1] \cup [e_1, e_2[, \quad 
 \Lambda_2 = ]e_1, e_2] \cup [e_2, e_0[, \quad 
 \Lambda_0 = ]e_2, e_0] \cup [e_0, e_2[ \] 
combine to a lift of a parametrization $[0,1] \to 
\partial\Delta^2 \cong \bS^1$ starting in the base point to a 
curve in $\hat K$. This implies that 
\[ \partial_1^{\hat K}([f]) 
= (\hat f_1^{-1} \hat f_0)\res_{]e_0, e_1[} 
= (\hat f_0\hat f_1^{-1} )\res_{]e_0, e_1[} = S_1([f^*\hat K]),\] 
where the next to last equality follows from the fact that 
$D$ is central in $\hat K$. 
This completes the proof. 
\end{prf}

\begin{cor}\mlabel{cor:lower-degree-analog}
If $Z$ is a $K(\Gamma,1)$-group 
and $Q \to \bS^2$ a $Z$-bundle, then the 
Chern class $\delta_1([Q]) \in \check H^2(\bS^2,\Gamma)$ 
satisfies 
\[ S_2(\delta_1([Q])) = -h([Q]) \in \pi_1(Z) \cong \Gamma.\] 
\end{cor}

\begin{prf}
 Applying the preceding proposition with $K=Z$ and $\hat{K}=\tilde Z 
\cong E\Gamma$ 
 shows the claim since $\partial_1^Z$ implements the isomorphism
 $\pi_{1}(Z)\cong \Gamma$.
\end{prf}

The following corollary now implies Theorem~\ref{thm:main}(b). 

\begin{cor} \mlabel{cor:5.3} Assume that $K$ is connected and that 
$\hat K$ is a central extension of $K$ by the discrete group~$D$. 
For a $K$-bundle $P \to X$ we then have 
\[ \alpha_2(\delta_1([P])) 
= -\partial_1^{\hat K} \circ \partial_2^P \: \pi_2(X) \to D.\] 
\end{cor}

\begin{prf} With Proposition~\ref{prop:lower-degree-analog},
the naturality of $\delta_1$ and Remark~\ref{rem:3.1} we obtain 
\begin{align*}
\alpha_2(\delta_1([P]))([\sigma]) 
&= S_2(\sigma^*\delta_1([P])) 
= S_2(\delta_1([\sigma^*P])) \\
&= -\partial_1^{\hat K} h(\sigma^*P) 
= -\partial_1^{\hat K}\circ \partial_2^P([\sigma]).
\qedhere  
\end{align*}
\end{prf}

\begin{rem} \mlabel{rem:1conn} (a) If $X$ is $1$-connected, then 
$H^2(X,D) \cong \Hom(\pi_2(X),D)$ by the Universal Coefficient Theorem 
and the Hurewicz Isomorphism $h_2 \: \pi_2(X) \to H_2(X)$. 
Therefore Corollary~\ref{cor:5.3} determines in this case 
the obstruction class completely. 

(b) If $X$ is not $1$-connected and $q_X \: \tilde X \to X$ is a 
simply connected  covering, then 
$\delta_1([q_X^*P]) =   q_X^*\delta_1([P])$ 
by Lemma~\ref{lem:0.2}, and Corollary~\ref{cor:5.3} thus leads with 
Lemma~\ref{lem:0} to 
\[ \alpha_2(\obs_{q_X^*P}([\hat K])) 
= -\partial_1^{\hat K} \circ \partial_2^P \circ \pi_2(q_X).\] 
As $\pi_2(q_X) \: \pi_2(\tilde X) \to \pi_2(X)$ is an isomorphism, 
(a) shows that Corollary~\ref{cor:5.3} determines the obstruction class of the 
$K$-bundle $q_X^*P$ over $\tilde X$. 

(c) 
From (b) we derive in particular that 
$\obs_{q_X^*P}([\hat K])$ vanishes if and only if 
$\partial_1^{\hat K} \circ \partial_2^P$ vanishes. For $\hat K = \tilde K$ 
the homomorphism $\partial_1^{\hat K}$ also is an isomorphism, so that 
$q_X^*P$ lifts to a $\tilde K$-bundle if and only if $\partial_2^P$ vanishes. 
\end{rem}

\begin{prop} \mlabel{prop:8.4} 
If $Z$ is $K(\Gamma,1)$-group, then the Chern class 
\[ \Ch \: \Bun(X,Z)\to \check H^2(X,\Gamma) \] 
defines an isomorphism of abelian groups mapping 
\[ \Bun(X,Z)_0 := \{ [P] \in \Bun(X,Z) \: \partial_2^P = 0\} \] 
onto the subgroup $\Lambda^2(X,\Gamma)$. 
If $q_X \: \tilde X \to X$ is the universal covering map, then 
$\partial_2^P = 0$ is equivalent to the triviality of the bundle 
$q_X^*P$ on $\tilde X$. 
\end{prop}

\begin{prf} We know already that $\Ch$ is an isomorphism of abelian 
groups. From $\Ch(P) = \obs^{\tilde Z}([P])$ and Remark~\ref{rem:1conn}(b) 
we further derive that $\alpha_2(\Ch(P))$ vanishes if and only if 
$\partial_2^P = 0$. Finally, Remark~\ref{rem:1conn}(c) 
shows that this in turn is equivalent to the triviality 
of the bundle $q_X^*P$ on~$\tilde X$.
\end{prf}

\begin{rem}\label{rem:5.6} Recall that for each topological space $X$ the elements in the group 
 $H_2(X)$ are represented by
maps $\Sigma \to X$, where $\Sigma$ is an orientable 
surface \cite[Ex.\ 1.4]{Cal09}.

Let $q \: \Sigma \to \bS^2$ be a map inducing an isomorphism 
$H_2(q) \: H_2(\Sigma) \to H_2(\bS^2)$. If 
$\Sigma$ is of genus $g$ and obtained from gluing the sides of a 
$4g$-gon, then such a map can be obtained by collapsing the whole 
boundary of this polygon to a point, which leads to a $2$-sphere. 
Let $\tilde S_2 \: \check H^2(\Sigma,\Gamma) \to \Gamma$ be the unique 
isomorphism with $\tilde S_2 \circ q^* = S_2$, so that 
we can evaluate $c \in \check H^2(X,\Gamma)$ 
on the homology class defined by a continuous map 
$\sigma \: \Sigma \to X$ by 
$\tilde S_2(\sigma^*c) \in \Gamma$. 

For a connected Lie group $K$, the classifying space $BK$ is 
$1$-connected, so that \cite[Cor.~13.16]{Br97a} 
implies that 
\[ \Bun(\Sigma, K) \cong [\Sigma,BK] \to 
H^2(\Sigma,\pi_2(BK)) \cong H^2(\Sigma,\pi_1(K)), 
\quad f \mapsto f^*[u] \]
is bijective, where 
\[ u \in H^2(BK,\pi_2(BK)) 
\cong \End(\pi_2(BK)) \cong \End(\pi_1(K)) \]
is the element corresponding to 
$\id_{\pi_1(K)}$.\begin{footnote}{Here we assume that 
$BK$ is homotopy equivalent to a a CW complex. This is the case for each metrizable
Lie group, cf.\ \cite[Lemma 4.4]{NSW11}}  
\end{footnote}
Since $H_1(\Sigma)$ is free and $H_2(\Sigma) \cong \Z$, the 
Universal Coefficient Theorem implies that 
\[ \Bun(\Sigma,K) \cong H^2(\Sigma,\pi_1(K)) \cong 
\Hom(H_2(\Sigma), \pi_1(K)) \cong \pi_1(K).\] 
The preceding discussion now implies that the pullback map  
\[ q^* \: \Bun(\bS^2,K) \to \Bun(\Sigma,K) \] 
is a bijection. 
If $h(P) \in \pi_1(K)$ is the characteristic class of the 
$K$-bundle $P$ over $\Sigma$ and $P \cong  q^*Q$ for a $K$-bundle $Q$ on $\bS^2$, then 
\[ \partial_1^{\hat K}h(P) = \partial_1^{\hat K}h(Q) =  
- S_2(\delta_1([Q]))
= - \tilde S_2(\delta_1([P])).\] 
For a general $K$-bundle $P \to X$ and the corresponding 
classifying map $f \: X \to BK$, we thus obtain for a 
continuous map $\sigma \: \Sigma \to X$: 
\begin{align*}
\tilde S_2(\sigma^*\delta_1([P])) 
&= \tilde S_2(\delta_1([\sigma^*P])) 
= -\partial_1^{\hat K} h(\sigma^*P)\\
&= -\partial_1^{\hat K}\circ \partial_2^{EK}[f \circ \sigma] 
= -\partial_1^{\hat K} \circ H_2(f)[\sigma].
\end{align*}
This means that the homomorphism $H_2(X) \to \pi_1(K)$ defined by 
the obstruction class $\delta_1([P]) \in H^2(X,\Gamma)$  coincides 
for $P \sim f^*EK$ with $-\partial_1^{\hat K} \circ H_2(f)$. 
\end{rem}

\subsection*{A degree zero analog}  

As a byproduct of Proposition~\ref{prop:lower-degree-analog}, 
we obtain the following analog of 
Theorem~\ref{thm:main}, resp., 
\eqref{eq:identity},  for the connecting map 
$\delta_0$.  For a based map 
$f \in C(X,K) = \check H^0(X,\uline K)$, we think of the 
maps 
\[ \pi_n(f) \: \pi_n(X) \to \pi_n(K)\] 
as  analogs of the connecting 
maps $\partial^P_n \: \pi_n(X) \to \pi_{n-1}(K)$ 
associated to $[P] \in \check H^1(X,\uline K)$ 
and recall that $\delta_0([f]) \in \check H^1(X,\uline Z) 
\cong \check H^2(X, \Gamma)$ corresponds to the $Z$-bundle 
$f^*\hat K$ over~$X$. 

\begin{prop}\mlabel{prop:lowdeg} 
For $f \in C(X,K)$ and the central extension 
$\hat K$ of $K$ by the $K(\Gamma,1)$-group $Z$, we have 
 \begin{equation}
   \label{eq:prob2}
\partial^{\hat K}_2 \circ \pi_2(f) = -\alpha_2(\iota_{1}(\delta_0(f))) \: 
\pi_2(X) \to \Gamma.  \end{equation}
\end{prop}

\begin{prf} Corollary~\ref{cor:lower-degree-analog}
implies that the negative of 
\begin{equation*}
 \alpha_{2}(\iota_{1}(\delta_{0}(f)))([\sigma])=
\alpha_{2}(\iota_{1}([f^{*}\hat{K}]))([\sigma])=
S_{2}(\iota_{1}([\sigma^{*}f^{*}\hat{K}]))=
S_{2}(\delta_{1}([\sigma^{*}f^{*}\hat{K}]))
\end{equation*}
(note that $\iota_{1}=\delta_{1}$ in this case)
%
equals 
\[ h([\sigma^*f^*\hat K]) 
= \partial_2^{\hat K}([f \circ \sigma]) 
= \partial_2^{\hat K}(\pi_2(f)([\sigma])).\] 
This proves \eqref{eq:prob2}. 
\end{prf}

\section{Proof of Theorem~\ref{thm:main}} 
\mlabel{sec:redux}

The first part in the proof of Theorem~\ref{thm:main} is to
reduce the problem to the case $X = \bS^3$. 
Our goal is to show the following equality of two group homomorphisms 
\begin{equation}
  \label{eq:identity}
\partial^{\hat K}_2 \circ \partial^P_3 = \alpha_3(\iota_{2}( \delta_1([P]))) \: 
\pi_3(X) \to \Gamma.
\end{equation}

\begin{lem} \mlabel{lem:redux} If \eqref{eq:identity} holds for 
$X=\bS^3$, then it holds for arbitrary~$X$. 
\end{lem}

\begin{prf} For a continuous map $\sigma \: \bS^3 \to X$, we obtain with the 
naturality of $\delta_1$ and $\iota_{2}$: 
\begin{align*}
\alpha_3(\iota_{2}^{X}(\delta_1^X([P])))([\sigma]) 
&= S_3(\sigma^*\iota_{2}^{X}(\delta_1^X([P])))\\
&= S_3\big(\iota_{2}^{\bS^{3}}(\delta_1^{\bS^3}([\sigma^*P]))\big) 
= \alpha_3\big(\iota_{2}^{\bS^{3}}(\delta_1^{\bS^3}([\sigma^*P]))\big)([\id_{\bS^3}]).
\end{align*}
If \eqref{eq:identity} holds for $\bS^3$, then 
\[ \partial_2^{\hat K} \circ \partial_3^{\sigma^*P} = 
\alpha_3(\iota_{2}^{\bS^{3}}(\delta_1^{\bS^3}([\sigma^*P]))) \: \pi_3(\bS^3) \cong \Z \to \Gamma,\] 
and  we obtain with Lemma~\ref{lem:0}
\begin{align*}
\alpha_3(\iota_{2}^{X}(\delta_1^X([P])))([\sigma]) 
&= \alpha_3\big(\iota_{2}^{\bS^{3}}(\delta_1^{\bS^3}([\sigma^*P]))\big)([\id_{\bS^3}])\\
&=\partial_2^{\hat K} \circ \partial_3^{\sigma^*P}([\id_{\bS^3}])
=\partial_2^{\hat K} \circ \partial_3^{P}([\sigma]),
\end{align*}
which implies \eqref{eq:identity}. 
\end{prf}

In view of Lemma~\ref{lem:redux}, 
it remains to show \eqref{eq:identity} for $X = \bS^3$. 
Since $\pi_3(\bS^3)$ is cyclic generated by the homotopy class of 
$\id_{\bS^3}$, we have to verify that for every $K$-bundle $P$ 
over $\bS^3$
\[ \partial_2^{\hat K} \circ \partial_3^{P}([\id_{\bS^3}]) 
= \partial_2^{\hat K}(h(P)) \] 
coincides with
\[S_3(\iota_{2}(\delta_1([P])))
= \alpha_3(\iota_{2}(\delta_1([P])))([\id_{\bS^3}]).\] 

We address this problem by verifying 
it for the bundles $P_f$ associated to a 
continuous map $f \: \bS^2 \to K$ in such a way that 
$h([P_f]) = [f]$. Since each $K$-bundle on $\bS^3$ arises from this 
construction, this will prove our theorem. 
In view of 
\[ \partial^{\hat K}_2([f]) 
= \partial^{f^*\hat K}_2([\id_{\bS^2}])
= h([f^*\hat K]) 
= - S_2(\delta_{1}([f^*\hat K]))= - S_2(\iota_{1}([f^*\hat K]))
,\] 
(Corollary~\ref{cor:lower-degree-analog}) and $[f^{*}\hat{K}]=\delta_{0}^{\bS^{2}}(f)$
(Proposition~\ref{Quadratdiagramm}) we have to verify that 
\begin{equation}
  \label{eq:sphercase}
-S_2(\iota _{1}(\delta_{0}^{\bS^{2}}(f))) = S_3(\iota_{2}(\delta_1^{\bS^3}([P_f]))).
\end{equation}
From Propositions \ref{Dreiecksdiagramm} we know that 
\[ S_3 \circ \susp{2} = S_2, \quad 
\susp{2} \circ \iota_1 = \iota_2 \circ \susp{1}\quad \mbox{ and } 
\quad \susp{0}(f) = [P_f].\] 
With Proposition~\ref{Quadratdiagramm}, we thus obtain 
\[\iota_2 \delta_1^{\bS^3}([P_f])
= \iota_2 \delta_1^{\bS^3}\susp{0}(f) 
=  -\iota_2(\susp{1}(\delta_0^{\bS^2}(f))
=  -\susp{2}(\iota_1(\delta_0^{\bS^2}(f))).\] 
Applying $S_3$ now completes the proof of 
Theorem~\ref{thm:main}.

\section{Aspherical obstruction classes} 
\mlabel{sec:6}

In this section we study those obstruction classes 
which cannot be detected by the corresponding homomorphisms 
on $\pi_2(X)$, resp., $\pi_3(X)$. These are the classes contained 
in the subgroups $\Lambda^k(X,\Gamma) = \ker \alpha_k \subeq 
\check H^k(X,\Gamma)$ 
of those cohomology classes inducing 
the trivial homomorphism $\pi_k(X) \to \Gamma$ ($k = 2,3$).
The Universal Coefficient 
Theorem (cf.\ \cite{Br97a}) leads to the short 
exact sequence 
\[ 0\to{\Ext}(H_{k-1}(X),\Gamma)\to H^k(X,\Gamma)\ssmapright{\beta_k}
 \Hom(H_k(X),\Gamma)\to 0\]
and the map $\beta_k$ satisfies 
\[ \alpha_k(c) = \beta_k(c) \circ h_k \: \pi_k(X) \to \Gamma \quad \mbox{ for } \quad 
c \in H^k(X,\Gamma).\] 
The image $\Sigma_k(X) := h_k(\pi_k(X))$ 
is called the {\it subgroup of spherical cycles} 
in $H_k(X)$.  Note that cohomology classes coming from elements of 
$\Ext(H_{k-1}(X),\Gamma)$ are in particular contained in 
$\Lambda^k(X,\Gamma)$ and that this inclusion may be proper. 
For $k = 2$ Hopf's Theorem asserts that 
$\Lambda^2(X,\Gamma) \cong H^2_{\rm grp}(\pi_1(X),\Gamma)$. 
For an interpretation of this result in terms of obstruction 
classes see Theorem~\ref{thm:7.7} below. 

\subsection{Central extensions by discrete groups}

If $K$ is connected, then $\partial_1^P$ is always trivial. 
Therefore we first  take a closer look at the set 
\[ \Bun(X,K)_0 := \{ [P] \in \Bun(X,K) \: \partial_2^P = 0\}\]  
because the vanishing of $\partial^P_1$ and $\partial_2^P$ implies 
that the corresponding obstruction classes are aspherical 
(cf.\ Theorem~\ref{thm:main}).
It turns out that, for every $[P] \in \Bun(X,K)_0$, the fundamental 
group $\pi_1(P)$ is a central extension of $\pi_1(X)$ by $\pi_1(K)$, 
and that the class $[\pi_1(P)] \in H^2_{\rm grp}(\pi_1(X),\pi_1(K))$ coincides, 
up to sign, with the lifting obstruction $\obs^{\tilde K}([P])$. 
This observation complements Theorem~\ref{thm:main} 
and leads to a description of all obstruction classes 
associated to flat central extensions of $K$ by $K(\Gamma,1)$-groups 
as  elements of 
$\Ext(H_2(X),\Gamma) \subeq H^3(X,\Gamma)$ (Theorem~\ref{thm:7.11}). 

Suppose that $X$ is connected and let $q_X \: \tilde X \to X$ denote a 
universal covering. We identify $\pi_1(X)$ with the group of deck 
transformations of $X$ and view $\tilde X$ as a $\pi_1(X)$-principal 
bundle. In this section $K$ denotes a connected locally contractible topological group. 

\begin{rem} \mlabel{rem:aspherspace} 
(a) If $\pi_2(X)$ or $\pi_1(K)$ vanishes, then 
$\Bun(X,K) = \Bun(X,K)_0$. 
This is in particular the case if $X$ is  {\it aspherical}, 
i.e., if $\pi_{k}(X)$ vanishes for $k\geq 2$.

Typical examples of aspherical spaces are: 
connected surfaces of positive genus, tori, solvmanifolds 
(the covering is a simply connected 
solvable Lie group) and  hyperbolic manifolds. 

(b) Suppose that $q_X^*P$ is trivial. 
Since the homomorphisms $\pi_n(q_X) \: \pi_n(\tilde X) \to \pi_n(X)$ 
are isomorphisms for $n > 1$, all the connecting maps 
$\partial_n^{P} \: \pi_n(X) \to \pi_{n-1}(K)$ vanish. 
For $n = 1$ we have a homomorphism  
$\partial_1^{P} \: \pi_1(X) \to \pi_0(K)$ which is also trivial 
if $K$ is connected. For a connected group $K$, we thus obtain 
$K$-bundles $P$ for which all connecting maps $\partial_d^P$ vanish. 
In view of Theorem~\ref{thm:main}, 
the corresponding obstruction classes are aspherical. 
\end{rem}

The following lemmas show that the obstruction classes defined by 
central extensions of $\pi_1(X)$ are aspherical. 

\begin{lem} \label{lem:obinlambda} 
If $K$ is discrete and $P \to X$ a principal $K$-bundle, then 
\[ \obs_P(\Ext_c(K,Z)) \subeq \Lambda^2(X,Z) \] 
for every abelian group $Z$. In particular, 
\[ \obs_{\tilde X}(H^2_{\rm grp}(\pi_1(X),Z)) \subeq \Lambda^2(X,Z).\] 
\end{lem}

\begin{prf} We have to show that $\sigma^*\obs_P([\hat K])$ 
vanishes for every continuous map $\sigma \: \bS^2 \to X$. 
In view of the naturality of the obstruction class, 
\[ \sigma^*\obs_P([\hat K]) 
= \obs_{\sigma^*P}([\hat K]), \]
and since the $K$-bundle 
$\sigma^*P \to \bS^2$ is trivial because 
$K$ is discrete and $\bS^2$ is simply connected, 
all corresponding obstruction classes vanish. 
\end{prf}

We now turn the universal covering space of a 
$K$-principal bundle into a principal bundle for a suitable extension 
of $K$ by the discrete group $\pi_1(P)$. 

\begin{lem} \mlabel{lem:7.1} Let $q \: P \to X$ be a connected $K$-bundle 
and $\alpha \: \pi_1(K) \to \pi_1(P)$ be the homomorphism defined by the choice 
of a base point in $P$. Then 
\[ K^\sharp := (\tilde K \times \pi_1(P))/\{(d,\alpha(d)^{-1}) \: d \in \pi_1(K)\} \] 
is an extension of $K$ by the discrete group $\pi_1(P)$ 
with $\pi_0(K^\sharp) \cong \pi_1(P)/\im(\alpha) \cong \pi_1(X)$ 
and the universal covering space $\tilde P$ of $P$ carries a natural 
$K^\sharp$-bundle structure such that $P \sim \tilde P/\pi_1(P)$ 
is an equivalence of $K$-bundles. 
\end{lem}

\begin{prf}  From basic covering theory we know that 
the $K$-action on $P$ lifts to a $\tilde K$-action on the universal covering 
space $\tilde P$. Here the action of $\pi_1(K) \cong \ker q_K$ is given by 
the homomorphism $\alpha \: \pi_1(K) \to \pi_1(P)$ obtained from the orbit map of a 
base point $p_0 \in P$, where we identify 
$\pi_1(P) \subeq \Homeo(\tilde P)$ with  the group of deck transformations of the 
covering~$q_P \: \tilde P \to P$. 
In particular, $D_0 := \ker\alpha$ acts trivially on $\tilde P$ 
and the group $\tilde K/D_0$ acts faithfully on $\tilde P$. 

Since $\tilde K$ is connected and normalizes the subgroup $\pi_1(P)\subeq \Homeo(\tilde P)$ 
whose orbits are discrete, it actually centralizes $\pi_1(P)$. 
We thus obtain an action of the group $K^\sharp$ because 
the discrete subgroup $\{(d,\alpha(d)^{-1}) \: d \in \pi_1(K)\}$ 
acts trivially. 

If $U \subeq X$ is an open $1$-connected subset for which we have a 
continuous section $\sigma \: U \to P$, then $\sigma$ also lifts to a 
continuous map $\tilde\sigma \: U \to \tilde P$. Now the action map 
leads to a local homeomorphism 
\[ U \times K^\sharp \to q_P^{-1}(U) \subeq \tilde P, \quad (x,k) \mapsto \tilde\sigma(x)k, \]
and since $\pi_1(P)$ acts freely on $\tilde P$, this actually is a 
$K^\sharp$-equivariant homeomorphism. Therefore $\tilde P$ is a $K^\sharp$-principal 
bundle and $q_P \: \tilde P \to P$ induces an equivalence 
$\tilde P/\pi_1(P) \to P$ of $K$-bundles. 
\end{prf}

\begin{rem} \mlabel{rem:7.2} (a) From the long exact homotopy sequence of the $K$-bundle 
$P$ we obtain the exact sequence 
  \begin{equation}
    \label{eq:longexact}
\pi_2(X) 
\ssmapright{\partial_2^P} \pi_1(K) 
\ssmapright{\alpha} \pi_1(P) 
\to \pi_1(X) \to \pi_0(K) = \{\1\},
  \end{equation}
so that $\partial_2^P = 0$ is equivalent to $\alpha$ being injective, resp., 
to $\tilde P$ carrying the structure of a $\tilde K$-principal bundle. 

(b) Since the $\tilde K$-action on $\tilde P$ commutes with $\pi_1(P)$, the subgroup 
$\alpha(\pi_1(K)) \subeq \pi_1(P)$ is central in $\pi_1(P)$. In particular, 
$\pi_1(P)$ is a central extension of $\pi_1(X)$ by $\im(\alpha)$, 
and if $\partial_2^P$ vanishes, then 
$[\pi_1(P)] \in H^2_{\rm grp}(\pi_1(X), \pi_1(K))$. 

(c) Let $K^+ := (K^\sharp)_0 \cong \tilde K/\im(\partial_2^P)$ 
be the identity component 
of~$K^\sharp$. 
Since $K^+$ is connected and $\tilde P$ is simply connected, 
the long exact homotopy sequence of the $K^+$-bundle $\tilde P \to 
\tilde P/K^+$ 
implies that this quotient is simply connected. It also 
is a $\pi_1(X)$-bundle and a covering of $X$, hence a model for the universal covering 
space $\tilde X$ of~$X$. 

(d) If the homomorphism $\pi_2(q) \: \pi_2(P) \to \pi_2(X)$ vanishes, then 
\eqref{eq:longexact} defines an exact four term sequence 
\begin{equation}
\label{eq:longexact2}
\{0\} \to \pi_2(X) 
\ssmapright{\partial_2^P} \pi_1(K) 
\ssmapright{\alpha} \pi_1(P) 
\to \pi_1(X) \to \{\1\}.
\end{equation}
This is the four terms sequence defined by the crossed module 
obtained from $\alpha$ and the trivial action of $\pi_1(P)$ on $\pi_1(K)$. 
In particular, the cohomology class of this crossed module leads to 
an element of $H^3_{\rm grp}(\pi_1(X), \pi_2(X))$ which can only be nonzero 
if $\partial_2^P \not=0$. 
%
\end{rem}

\begin{lem} \mlabel{lem:sumobs} Let $q_j \: P_j \to X$ be $K_j$-bundles for $j =1,2$. 
If $\hat K_j$, $j =1,2$, are central $Z$-extensions of $K_j$ and 
$\mu_Z \: Z \times Z \to Z$ is the multiplication map, then 
\[ Q := (\mu_Z)_*(\hat K_1 \times \hat K_2) \] 
is a central extension of $K_1 \times K_2$ by $Z$. 
For the obstruction class of the $K_1 \times K_2$-bundle 
$P_1 \times_X P_2$ over $X$ we then have 
\[ \obs_{P_1 \times_X P_2}([Q]) = \obs_{P_1}([\hat K_1]) + \obs_{P_2}([\hat K_2]) 
\in \check H^2(X,\uline Z).\] 
\end{lem}

\begin{prf} This is an immediate consequence of the definition of the obstruction 
class. 
\end{prf}

The following proposition is a key tool which relates certain 
obstruction classes for general connected group $K$ to obstruction 
classes for central extensions of $\pi_1(X)$.

\begin{prop} \mlabel{prop:8.3a} 
If $K$ is connected, $P \to X$ a $K$-bundle, $D := \im(\alpha) \subeq \pi_1(P)$ 
and $K^+ := \tilde K/\im(\partial_2^P)$, then 
\[\obs_{P}([K^+]) 
= -\obs_{\tilde X}([\pi_1(P)])  \in \Lambda^2(X,D). \] 
If $\partial_2^P = 0$, then $D= \pi_1(K)$ and we have in particular 
\[\obs_{P}([\tilde K]) = -\obs_{\tilde X}([\pi_1(P)]) \in \Lambda^2(X,\pi_1(K)).\] 
\end{prop}
  
\begin{prf} We apply Lemma~\ref{lem:sumobs} with 
$Z = D$, $K_1 = K$, $\hat K_1 = K^+$, 
$P_1 = P$, $K_2 = \pi_1(X)$, $\hat K_2 = \pi_1(P)$ and 
$P_2 = \tilde X$ (Remark~\ref{rem:7.2}(b)). Then $Q \cong K^\sharp$ 
as a central 
extension of $K \times \pi_1(X)$ by $D$ (Lemma~\ref{lem:7.1}). 

As $\tilde P/K^+ \sim \tilde X$ as $\pi_1(X)$-bundles and 
$\tilde P/\pi_1(P) \sim P$ as $K$-bundles, it follows that 
\[ \tilde P/D \sim P \times_X \tilde X = q_X^*P \] 
as $(K \times \pi_1(X))$-bundle. In particular, 
the $K$-bundle $q_X^*P$ lifts to a $K^+$-bundle, so that 
\[ 0 = \obs_{q_X^*P}([K^+]) = \obs_P([K^+]) + \obs_{\tilde X}([\pi_1(P)]), \] 
where the second equality follows from Lemma~\ref{lem:sumobs}. 
\end{prf}

We now have all the tools to give a complete description of the 
aspherical obstruction classes for the case where $Z = D$ is discrete. 
If $X$ is $1$-connected, then this has already been done in 
Corollary~\ref{cor:5.3} and the subsequent Remark~\ref{rem:1conn}. 

\begin{prop} \mlabel{prop:7.12} 
Suppose that $K$ is connected, $D$ is a discrete abelian group,  
 $\hat K$ is a central extension of $K$ by $D$ and that $P$ is a $K$-bundle 
over $X$. Then $\hat K\cong \gamma_*\tilde K$ for 
$\gamma = \partial_1^{\hat K}$ and 
$\delta_1([P])$ is aspherical if and only if 
$\im(\partial_2^P) \subeq \ker \gamma$. 
In this case the homomorphism $\gamma \: \pi_1(K) \to D$ factors through a 
homomorphism $\oline\gamma \: \coker(\partial_2^P) \to D$ and 
\[ \obs_P([\hat K]) 
= - \oline\gamma_*\obs_{\tilde X}([\pi_1(P)]) \in \Lambda^2(X,D),\] 
where $\pi_1(P)$ is considered as a central extension of $\pi_1(X)$ 
by $\im(\partial_2^P)$. 
\end{prop}

\begin{prf} From Remark~\ref{rem:3.1a}(a) we know that 
$\hat K \cong \gamma_*\tilde K$ for some 
$\gamma = \partial_1^{\hat K} 
\in \Hom(\pi_1(K),D)$. 
Since 
\[ \alpha_2(\delta_1([P])) 
= - \partial_1^{\hat K} \circ \partial_2^P
= - \gamma \circ \partial_2^P\] 
 by  Theorem~\ref{thm:main}(b), 
the obstruction class $\delta_1([P])$ is 
aspherical if and only if $\im(\partial_2^P) \subeq \ker\gamma$. 
In this case we obtain $\oline\gamma$ by factorization of $\gamma$, and 
for $K^+ := \tilde K/\im(\partial_2^P)$ we obtain
$\hat K \cong \oline\gamma_* K^+$. 
Therefore we can use 
Proposition~\ref{prop:8.3a} to obtain 
\[ \obs_P([\hat K]) 
=  \obs_P([\oline\gamma_*K^+]) 
=  \oline\gamma_*\obs_P([K^+]) 
= - \oline\gamma_*\obs_{\tilde X}([\pi_1(P)]).\qedhere\] 
\end{prf}

\begin{rem} Keeping the $K$-bundle $P \to X$ fixed, the preceding 
proposition can be used as follows to decide if two central 
$D$-extensions $\hat K_1$ and $\hat K_2$ have the same obstruction 
class. A necessary condition is that 
\[ \partial_1^{\hat K_1} \circ \partial_2^P 
= -\alpha_2(\obs_P([\hat K_1]))
= -\alpha_2(\obs_P([\hat K_2]))
= \partial_1^{\hat K_2} \circ \partial_2^P.\]
If this condition is satisfied, then 
the central extension $[\hat K_3] := [\hat K_1] - [\hat K_2]$ 
has an aspherical extension class which can be computed with 
Proposition~\ref{prop:7.12}. 
\end{rem}

The following theorem provides an obstruction theoretic version 
of Hopf's Theorem (cf.\ Example~\ref{ex:2.6}(c)).

\begin{theo} \mlabel{thm:7.7} For every abelian group $D$, 
\[ \obs_{\tilde X} \: H^2_{\rm grp}(\pi_1(X),D) \to \Lambda^2(X,D) \] 
is an isomorphism of abelian groups. 
\end{theo}

\begin{prf} From Lemma~\ref{lem:obinlambda} we  know that 
$\im(\obs_{\tilde X})) \subeq \Lambda^2(X,D)$. 

To see that $\obs_{\tilde X}$ is injective, suppose that 
$\obs_{\tilde X}([\hat\pi_1(X)])= 0$, i.e., that the 
$\pi_1(X)$-bundle $\tilde X$ lifts to a 
$\hat\pi_1(X)$-bundle $\hat q \: \hat X \to X$ 
with $\hat X/D \cong \tilde X$. 
Since $\hat q$ is a covering of $X$, it is 
associated to $\tilde X$ by a homomorphism 
$\gamma \: \pi_1(X) \to \hat\pi_1(X)$, and this homomorphism 
splits the extension $\hat\pi_1(X)$ of $\pi_1(X)$. Hence 
$\obs_{\tilde X}$ is injective. 

Now we show that $\obs_{\tilde X}$ is surjective. 
Let $Z$ be a $K(D,1)$-group. 
In view of Proposition~\ref{prop:8.4}, it suffices to show that 
every Chern class $\Ch(P)$, $[P] \in \Bun(X,Z)_0$, 
is contained in the range of $\obs_{\tilde X}$. 
In view of $\Ch(P) = \obs_{P}([\tilde Z])$ (Example~\ref{ex:2.6}(a)), 
this follows from Proposition~\ref{prop:8.3a}. 
\end{prf}

\begin{cor} \mlabel{cor:8.5a} 
If $K$ is connected, then we have a sequence of maps  
\begin{equation}
  \label{eq:tiltriv}
\Bun(X,\tilde K)_0 \to \Bun(X,K)_0 
\sssmapright{\zeta} H^2_{\rm grp}(\pi_1(X), \pi_1(K)), \quad 
\zeta([P]) = [\pi_1(P)], 
\end{equation}
which is exact in the sense that 
$\zeta([P]) = 0$ is equivalent to existence of a lift 
of $P$ to a $\tilde K$-bundle. 
\end{cor}

\begin{prf} From Proposition~\ref{prop:8.3a} we know that 
$\obs_P([\tilde K]) = - \obs_{\tilde X}([\pi_1(P)]),$  
so that the assertion follows from Theorem~\ref{thm:7.7}. 
\end{prf}

For abelian groups, the preceding observation can be refined as follows. 

\begin{prop} \mlabel{prop:8.7} 
For a  $K(\Gamma,1)$-group $Z$,  the following assertions hold: 
\begin{description}
\item[\rm(i)]  Any $[P] \in \Bun(X,Z)_0$ satisfies 
$\Ch(P) =  -\obs_{\tilde X}([\pi_1(P)]).$
\item[\rm(ii)]  
$\zeta \: \Bun(X,Z)_0 \to H^2_{\rm grp}(\pi_1(X), \Gamma), \zeta([P]) = [\pi_1(P)]$ 
is an isomorphism of abelian groups. 
\end{description}
\end{prop}

\begin{prf} (i) follows from $\Ch(P) = \obs_{P}([\tilde Z])$ (Example~\ref{ex:2.6}(a)) 
and Proposition~\ref{prop:8.3a}.

(ii) From (i) we derive that 
\[ \obs_{\tilde X} \circ \zeta = - \Ch \: \Bun(X,Z)_0 \to \Lambda^2(X,D).\] 
Since $\Ch$ and $\obs_{\tilde X}$ are isomorphisms by 
Proposition~\ref{prop:8.4} and Theorem~\ref{thm:7.7}, $\zeta$ also is an isomorphism.
\end{prf}

\subsection{Flat central extensions} 
\mlabel{subsec:flat} 

After the detailed discussion of central extensions 
by discrete groups, we now turn to the larger class of flat central 
extensions. Here the case where $Z$ is a $K(\Gamma,1)$-group 
leads to aspherical obstruction classes in $\Lambda^3(X,\Gamma)$. 
We start with the case of bundles $P$ for which $\partial_2^P$ vanishes 
and discuss general bundles over $1$-connected spaces in the following 
subsection.

If $\hat K = \gamma_* \tilde K$ is a flat central extension of $K$ by $Z$
and $[P] \in \Bun(X,K)_0$, 
then we can use Proposition~\ref{prop:7.12} and the naturality of the 
obstruction class (Lemma~\ref{lem:0.2}(c)) to obtain the relation 
\[ \obs_P([\hat K]) 
= \gamma_*\obs_P([\tilde  K]) 
= -\gamma_* \obs_{\tilde X}([\pi_1(P)]).\] 
To evaluate this formula, we first recall that 
$\obs_{\tilde X}([\pi_1(P)]) \in \check H^2(X,\pi_1(K))$ and that 
$\gamma_*$ denotes the natural map 
\[ \gamma_* \:  \check H^2(X,\pi_1(K)) 
\to \check H^2(X,\uline{Z}) \cong \check H^3(X,\Gamma).\] 
Writing $\iota^Z = \id_Z \: Z_d \to Z$ for the continuous bijection, where 
$Z_d$ denotes the discrete group $Z$, then $\gamma$ can be factorized as 
$\gamma = \iota^Z \circ \gamma_d$, which accordingly leads to 
$\gamma_* = \iota^Z_* \circ (\gamma_d)_*$ with 
\[ (\gamma_d)_* \: H^2(X,\pi_1(K)) \to H^2(X,Z) \quad \mbox{ and } \quad 
\iota^Z_* \: H^2(X,Z) \to H^3(X,\Gamma).\] 

If $Z$ is a $K(\Gamma,1)$-group, then the functors 
$\Ext(\cdot,Z)$ and $\Ext(\cdot,\tilde Z)$ vanish, so that the 
Universal Coefficient Theorem implies that 
\[ H^2(X,Z) \cong \Hom(H_2(X),Z) \quad \mbox{ and } \quad 
 H^2(X,\tilde Z) \cong \Hom(H_2(X),\tilde Z).\] 
Accordingly, the connecting map 
$\delta_2 \: H^2(X,Z) \to H^3(X,\Gamma)$ 
leads to the injection 
\[ \Ext(H_2(X),\Gamma) \cong \im(\delta_2) \into H^3(X,\Gamma) \] 
from the Universal Coefficient Theorem. This implies that 
\[ \im(\iota^Z_*) = \Ext(H_2(X),\Gamma) \subeq H^3(X,\Gamma).\] 
We thus arrive that the following theorem complementing 
Theorem~\ref{thm:main} by providing similar information 
on aspherical obstruction classes for $\partial_2^P = 0$.

\begin{theo} \mlabel{thm:7.11} 
If $K$ is connected, $[P] \in \Bun(X,K)_0$, $Z$ is a $K(\Gamma,1)$-group 
and $\hat K = \gamma_*\tilde K$ a flat extension of $K$ by $Z$, then the corresponding 
obstruction class is given by 
\[ \obs_P([\hat K]) 
= -\iota^Z_*\obs_{\tilde X}((\gamma_d)_*[\pi_1(P)]) \in\Ext(H_2(X),\Gamma) 
\subeq H^3(X,\Gamma).\] 
\end{theo}

\begin{rem}
If the abelian group $\pi_1(K)$ is finitely generated, which is the case 
for a finite dimensional Lie group~$K$, then 
$\pi_1(K) \cong \Z^d \oplus F$ for some $d \in \N_0$ and a finite group $F$. 
For $Z= \T$ and $\Gamma = \Z$ we thus obtain a homomorphism 
\[ \oline\xi \: 
\Ext(\pi_1(K),\Z) 
= \Ext(F,\Z) 
\cong \Hom(F,\T) \cong F \to H^3(X,\Gamma).\] 

Note that, for every finite  CW-complex $X$, the groups $H_k(X)$ are finitely 
generated, so that 
\[ H^3(X,\Z) \cong \Ext(H_2(X),\Z) \oplus \Hom(H_3(X),\Z), \] 
where $\Hom(H_3(X),\Z)$ is free and 
\begin{align*}
\Tor(H^3(X,\Z)) 
&\cong \Ext(H_2(X),\Z) 
\cong \Ext(\Tor(H_2(X)),\Z) \\
&\cong \Hom(\Tor(H_2(X)),\T) \cong \Tor(H_2(X)) 
\end{align*}
is finite. 
\end{rem}

\begin{rem} \mlabel{rem:7.14} (Flat bundles) 
(a) Let $\gamma \: \pi_1(X) \to K$ be a group homomorphism 
and 
\[ P^\gamma := \gamma_*\tilde X = (\tilde X \times K)/\pi_1(X) \] 
denote the 
corresponding flat bundle, which is the quotient of 
$\tilde X \times K$ by  the right action of $\pi_1(X)$ by 
$(x,k).d = (xd, \gamma(d)^{-1}k)$. 
Then $\pi_1(P^\gamma)$ is a central extension of 
$\pi_1(X)$ by $\pi_1(K)$, and since $\tilde{P^\gamma} \cong \tilde X \times \tilde K$, 
it follows that 
\[ \pi_1(P^\gamma) \cong 
\{ (d,\tilde k) \in \pi_1(X) \times \tilde K \: \gamma(d) = q_K(\tilde k)\} 
= \gamma^*\tilde K,\] 
as a central extension of $\pi_1(X)$ by $\pi_1(K)$. 
In view of Corollary~\ref{cor:8.5a}, the bundle 
$P$ lifts to a $\tilde K$-bundle over $X$ if and only if 
$\zeta([P^\gamma]) = [\pi_1(P^\gamma)]$ vanishes, which in 
turn is equivalent to the existence of a lift 
$\tilde\gamma \: \pi_1(X) \to \tilde K$. 
In particular, any flat $K$-bundle that lifts to some $\tilde K$-bundle 
lifts to a flat $\tilde K$-bundle. 

(b) Let $Z$ be a $K(\Gamma,1)$-group. 
If $P^\gamma \to X$ is a flat $Z$-bundle defined by a homomorphism 
$\gamma \: \pi_1(X) \to Z$, then the Hurewicz Theorem implies that 
$\gamma$ factors through a homomorphism 
\[\oline\gamma \: H_1(X) \cong \pi_1(X)/(\pi_1(X),\pi_1(X)) \to Z.\] 
Therefore 
\[ \zeta(P^\gamma) = [\pi_1(P^\gamma)] 
= \gamma^*\tilde Z \] 
lies in the subgroup $\Ext(H_1(X),\Gamma)\subeq H^2(X,\Gamma)$, 
 parametrizing those central 
extensions pulled back from abelian extensions of $H_1(X)$ by the 
quotient map $h_1 \: \pi_1(X) \to H_1(X)$. 
Since $\tilde Z$ is divisible, $\Ext(H_1(X),\tilde Z)$ vanishes, and therefore 
the exact Hom-Ext sequence for abelian groups shows that 
the connecting map 
\[ \delta_1 \: \Hom(H_1(X),Z) \to \Ext(H_1(X),\Gamma) \] 
is surjective, i.e., every class in 
$\Ext(H_1(X),\Gamma)$ can be represented by a homomorphism 
$\gamma \: \pi_1(X) \to Z$ as $[\oline\gamma^*\tilde Z]$. 
This implies that 
\[ \Ch(\Bun(X,Z)_{\rm flat}) = \Ext(H_1(X),\Gamma) \subeq 
\Ch(\Bun(X,Z)_0) = \Lambda^2(X,\Gamma).\] 
\end{rem}

\subsection{Aspherical classes for $1$-connected spaces} 

If $X$ is $1$-connected, then  $\Lambda^2(X,\Gamma) 
\cong H^2_{\rm grp}(\pi_1(X),\Gamma)$ vanishes. 
Therefore we only have non-zero aspherical obstruction classes 
if $Z$ is not discrete. Let us first recall the structure of 
$\Lambda^3(X,\Gamma)$. 

\begin{rem} \mlabel{rem:6.1} (a) If $\pi_2(X)$ vanishes, then it follows from 
\cite[Thm.~II]{EML45} that 
\[ H_3(X)/\Sigma_3(X) \cong H_3(\pi_1(X))\]  
is the third homology group of $\pi_1(X)$. 
 If, in addition, $X$ is $2$-connected, then 
$h_3 \: \pi_3(X) \to H_3(X)$ is an isomorphism and 
$H_2(X)$ vanishes. Therefore 
$\Lambda^3(X,\Gamma) = \{0\}$.

(b) If $X$ is $1$-connected, then 
\cite[Thm.\ II$^m$, p.~509]{EML45} shows that 
\[ H_3(X)/\Sigma_3(X) \cong H_3(K(\pi_2(X),2)) = \{0\},\]  
where the last equality is due to the fact that 
the homology algebra of any topological group of type 
$K(\pi_2(X),2)$ (which always exists by \cite{Mi67}) 
is generated by $\pi_2(X)$ in degree $2$, 
so that all odd degree classes vanish. We conclude that 
\begin{equation}
  \label{eq:e1}
\Lambda^3(X,\Gamma) \cong \Ext(H_2(X),\Gamma). 
\end{equation}
Since the Hurewicz homomorphism $h_2 \: \pi_2(X) \to H_2(X)$ 
is an isomorphism, we also have 
\begin{equation}
  \label{eq:e2}
\Lambda^3(X,\Gamma) \cong \Ext(\pi_2(X),\Gamma). 
\end{equation}
\end{rem}

The preceding discussion shows that aspherical obstruction 
classes for bundles over $1$-connected spaces are most naturally 
represented by central extensions of $\pi_2(X)$. Therefore the 
problem is to find a natural way to describe this extension 
in terms of $\hat K$ and $P$. From Theorem~\ref{thm:main}(a) 
we know that $\delta_1([P])$ is aspherical if and only if 
$\partial_2^{\hat K} \circ \partial_3^P = 0$, which is always 
the case if $K$ is a finite dimensional Lie group. 

\begin{prop} If $X$ is $1$-connected, $K$ is connected 
and $\hat K = \gamma_*\tilde K$ a flat central extension of $K$ by the 
$K(\Gamma,1)$-group $Z$, then 
\[ \obs_P([\hat K]) = - (\gamma_d\circ \partial_2^P)^* \tilde Z 
\in \Ext(\pi_2(X), \Gamma) \cong \Lambda^3(X,\Gamma),\] 
where $\gamma_d \: \pi_1(K) \to Z_d$ is $\gamma$, considered as a homomorphism 
into the discrete group $Z_d$ underlying~$Z$.
\end{prop}

\begin{prf} By assumption, $\hat K = \gamma_*\tilde K$ for 
$\gamma = \partial_1^{\hat K}\: \pi_1(K) \to D$ (Remark~\ref{rem:3.1a}). 
Since $\tilde K$ is a central extension by a discrete group, 
Theorem~\ref{thm:main}(b) tells us that 
\[ \alpha_2(\obs_P([\tilde K])) 
= - \partial_1^{\tilde K} \circ \partial_2^P
= - \partial_2^P\] 
because $\partial_1^{\tilde K} = \id_{\pi_1(K)}$. 
Since $X$ is $1$-connected, 
\[ H^2(X,\pi_1(K)) \cong 
\Hom(\pi_2(X),\pi_1(K)),\] 
so that we have a complete 
description of $\obs_P([\tilde K])$. 
Using  the factorization $\gamma = \iota^Z \circ \gamma_d$ 
(cf.\ Subsection~\ref{subsec:flat}), 
we now obtain 
\begin{align*}
\obs_P([\hat K]) 
&= \obs_P([(\iota^Z)_*(\gamma_d)_*\tilde K])
= (\iota^Z)_* (\gamma_d)_* \obs_P([\tilde K]), 
\end{align*}
where 
\begin{align*}
(\iota^Z)_* &\: \Hom(\pi_2(X), Z) \cong \Hom(H_2(X), Z) \cong H^2(X,Z) \\
&\to 
 \Ext(H_2(X),\Gamma) \subeq H^3(X,\Gamma), \quad 
\alpha \mapsto [\alpha^*\tilde Z] 
\end{align*}
is the natural homomorphism which can be identified with the 
connecting homomorphism 
$\delta_2 \: \Hom(H_2(X),Z) \to \Ext(H_2(X),\Gamma)$ 
in the exact Hom-Ext sequence obtained for $H_2(X)$ from the short 
exact sequence $\Gamma \into \tilde Z \onto Z$. 

Identifying $H^2(X,Z_d)$ with $\Hom(H_2(X),Z_d)$ and 
$\Lambda^2(X,Z_d)$ with \break $\Ext(H_2(X),Z_d)$, we 
thus obtain 
\begin{align*}
\obs_P([\hat K]) 
= - (\iota^Z)_* (\gamma_d)_* \partial_2^P 
= - (\iota^Z)_* (\gamma_d \circ \partial_2^P) 
= - (\gamma_d \circ \partial_2^P)^*\tilde Z,
\end{align*}
and this completes the proof. 
\end{prf}

\begin{rem} It is also instructive to take a closer look at the 
obstruction classes of the universal bundle 
$q \: EK \to BK$ of a connected locally contractible topological group~$K$. 
We assume that $BK$ is locally contractible, which is the case for 
each metrizable Lie group (\cite[Lemma 4.4]{NSW11}). 
As $BK$ is $1$-connected, $H_2(BK) \cong \pi_2(BK) \cong \pi_1(K)$, so that the 
Universal Coefficient Theorem leads to 
\[ H^2(BK,\Gamma) \cong \Hom(\pi_1(K),\Gamma). \] 
Moreover, $\pi_3(BK) \cong \pi_2(K)$ and 
$\Lambda^3(BK,\Gamma) \cong \Ext(H_2(BK),\Gamma) \cong \Ext(\pi_1(K),\Gamma)$ 
(Remark~\ref{rem:6.1}(b)). 
Therefore the characteristic data determining cohomology classes in 
degrees $2$ and $3$ consists of elements of 
\[ \Hom(\pi_1(K),\Gamma), \quad 
\Hom(\pi_2(K),\Gamma) \quad \mbox{ and } \quad \Ext(\pi_1(K),\Gamma).\] 

For a central extension of $K$ by a discrete group $D$ we obtain 
with Theorem~\ref{thm:main}(b) the obstruction class 
\[ \obs_{EK}([\hat K]) = -\partial_1^{\hat K} \in H^2(BK,D) \cong \Hom(\pi_1(K),D).\]  
For a flat central extension $\hat K = \gamma_*\tilde K$ by 
the $K(\Gamma,1)$-group $Z$, defined by 
the homomorphism $\gamma \: \pi_1(K) \to Z$, we have 
\[ \obs_{EK}([\gamma_*\tilde K]) = -\delta(\gamma) \in 
\Ext(\pi_1(K),\Gamma).\] 
For a general central extension of $K$ by $Z$, 
Theorem~\ref{thm:main} implies that 
\[ \alpha_3(\obs_{EK}([\hat K])) = \partial_2^{\hat K} \: \pi_2(K) \to \Gamma \] 
because $\partial_3^P = \id_{\pi_2(K)}$. 
Note that these objects describe precisely the data discussed in 
Remark~\ref{rem:3.4b}. 

The class $\obs_{EK}([\hat K])$ is aspherical 
if and only if $[\hat K] \in \Bun(K,Z)_0$. If this is the case, then 
we know from Remark~\ref{rem:0.4} that that 
$\obs_{EK}([\hat K]) = \obs_{EK}([\gamma_*\tilde K])$ for some 
homomorphism $\gamma \: \pi_1(K) \to Z$, and in this case the obstruction 
class is $-\delta(\gamma) \in \Ext(\pi_1(K),\Gamma)$. 

Since a general $K$-bundle $P \to X$ is equivalent to $f^*EK$ for some 
continuous map $f \: X \to BK$, we can obtain the corresponding obstruction 
classes as pullbacks of the classes of $BK$. 
\end{rem}

From Remark~\ref{rem:6.1}(b), Theorem~\ref{thm:main}  
and the Universal Coefficient Theorem, 
we obtain in particular the following variant of
\cite{MurStev}.

\begin{prop} \mlabel{prob:6.9} 
If $X$ is $1$-connected and $\partial^{\hat K}_2 =0$, 
then  $\obs_P([\hat{K}])$ 
 is contained in $\Ext(H_{2}(X),\Gamma)\cong \Lambda^3(X,\Gamma)$. 
If $H_2(X)$ and $H_3(X)$ are finitely generated and 
$\Gamma$ is free, then it is a torsion element. 
\end{prop}

If $H_2(X)$ and $H_3(X)$ are finitely generated, then 
\[ \Tor(H^3(X,\Z)) 
\cong \Ext(H_2(X),\Z) \cong \Ext(\Tor(H_2(X)),\Z) \] 
follows from the Universal Coefficient Theorem, and this group is finite. 
Its elements can be considered as lifting obstructions 
of homomorphisms $\gamma \: H_2(X) \to \T \cong \R/\Z$.    

\begin{rem} \mlabel{rem:7.21} 
For a $1$-connected space $X$, a 
natural way to represent  classes in $\Ext(H_2(X),\Gamma)$ 
as obstruction classes is to consider a 
$K(H_2(X),1)$-group $H$ and write 
$q \: \hat X \to X$ for the $H$-bundle with Chern class 
\[ \Ch(\hat X) = \id_{H_2(X)} \in \End(H_2(X)) \cong H^2(X,H_2(X)). \] 

(a) As $\Ch(\hat X) = \obs^{\tilde H}([\hat X])$, we then obtain for a 
$K(\Gamma,1)$-group $Z$, 
a homomorphism $\gamma \: H_2(X) \to Z$, and the factorization 
$\gamma = \iota^Z \circ \gamma_d$, the relation 
\begin{align*}
\obs_{\hat X}(\gamma_*\tilde H) 
&=  \obs_{\hat X}(\iota^Z_*(\gamma_d)_*\tilde H) 
=  \iota^Z_* (\gamma_d)_* \obs_{\hat X}(\tilde H)\\ 
&=  \iota^Z_* (\gamma_d)_* \id_{H_2(X)} 
=  \iota^Z_* \gamma_d 
= \gamma_d^*{\tilde Z} \in \Ext(H_2(X),\Gamma).\end{align*}
If $H_2(X)$ is free and finitely generated, 
then $H$ can be realized as a torus but then 
$\Ext(H_2(X),\Gamma)$ vanishes. Therefore the preceding 
construction produces no non-trivial obstruction 
classes if $H$ is an abelian Lie group. 

(b) (Geometric interpretation of $\partial_3^P$) 
As $\partial_2^{\hat X} \: \pi_2(X) \cong H_2(X) \to \pi_1(H)$ 
is an isomorphism and $\pi_2(H)$ vanishes, the long exact homotopy sequence 
of $\hat X$ implies that $\hat X$ is $2$-connected. 
Therefore $H^3(\hat X,\Gamma) \cong \Hom(\pi_3(X),\Gamma)$. 

Suppose that $\hat K$ is a central $Z$-extension of $K$ 
for which $\partial_2^{\hat K}$ is an isomorphism. 
For a $K$-bundle $P \to X$ Theorem~\ref{thm:main}(a) now shows that 
\[ \alpha_3(\delta_1([q^*P]))
= \alpha_3(\delta_1([P])) \circ \pi_3(q)
= \partial_2^{\hat K} \circ \partial_3^P \circ \pi_3(q),\] 
and since $\pi_3(q) \: \pi_3(\hat X) \to \pi_3(X)$ is an isomorphism, 
we see that $\delta_1([q^*P])$ vanishes if and only if 
$\partial_3^P = 0$. This means that the vanishing of 
$\partial_3^P$ is equivalent to the existence of a lift of the 
$K$-bundle $q^*P$ to a $\hat K$-bundle. 
\end{rem}

As we shall see below, one obtains non-trivial aspherical obstruction 
classes from non-abelian Lie groups $K$ with finite fundamental 
group. As 
\[ H^2(X,\pi_1(K)) \cong \Hom(\pi_2(X),\pi_1(K))\] for a 
$1$-connected space $X$, the main point is to find 
$K$-bundles $P$ for which $\partial_2^P$ is non-trivial.

\begin{rem} \mlabel{rem:7.19} 
(a) If 
$K$ is a connected Lie group for which $H^1(\fk,\R)$ and $H^2(\fk,\R)$ vanish, 
we have 
\[ \Hom(\pi_1(K),\T) \cong  \Ext_s(K,\T)\] 
(cf.\ \cite[Thm.~7.12]{Ne02}). For $K = \PU_n(\C) \cong \SU_n(\C)/C_n \1$ 
this implies that 
\[ \Ext_s(\PU_n(\C),\T) \cong \Hom(C_n,\T) \cong C_n \] 
is  cyclic of order $n$
and the central $\T$-extension $\hat K = \U_n(\C)$ corresponds to a generator. 
For a $K$-bundle $P \to X$, the corresponding obstruction class 
$\obs_P([\hat K])$ is therefore contained in 
$\Tor(H^3(X,\Z)) \cong \Ext(H_2(X),\Z)$. 
According to \cite{Gr64}, all torsion classes are of this form. Indeed,
the {\it Brauer group} consists of equivalence classes of Azumaya algebras.
Grothendieck identifies equivalence classes of Azumaya algebras and 
${\rm PGl}_n$-bundles which do not lift to ${\rm Gl}_n$. He cites then in 
\cite{Gr64} a theorem of Serre which shows that for a finite CW complex $X$, 
the Brauer group is isomorphic to $\Tor(H^3(X,\Z))$.

(b) Let $K=\PU(\cH)$ be the projective unitary group of
an infinite dimensional separable Hilbert space $\cH$ 
and $\hat K := \U(\cH)$, which is a non-trivial central $\T$-extension 
for which 
\[ \partial_2^{\hat K} \: \pi_2(\PU({\mathcal H})) \to 
\pi_1(\T) \cong \Z \] 
is an isomorphism. 

For any compact locally contractible space $X$, it is known 
that any element of $H^3(X,\Z)$ is the Dixmier--Douady invariant 
$\obs^{\U(\cH)}([P])$ 
of a principal $\PU({\mathcal H})$-bundle $P\to X$ 
(see \cite{Dix63} or \cite[Thm.~5.1]{Sch}).\begin{footnote}
{Since the classifying space $B\PU(\cH)$ is an 
Eilenberg--MacLane space of type $K(\Z,3)$, Huber's Theorem 
(\cite{Hu61}) also provides a bijection $\Bun(X,\PU(\cH)) 
\cong [X,B\PU(\cH)] \cong \check H^3(X,\Z)$ for any $k$-space $X$.}
\end{footnote}
The obstruction class 
$\obs^{\U(\cH)}([P])$ is  non-zero 
for each non-trivial $\PU(\cH)$-bundle over $X$. 
\end{rem}

\begin{ex} For $X = \T^3$ the group 
$\pi_3(X) \cong \pi_3(\tilde X) = \pi_3(\R^3)$ 
is trivial and $\Lambda^3(X,\Z) = 
H^3(X,\Z) \cong \Z$ because $X$ is an orientable 
$3$-manifold. Hence $\partial_3^P = 0 = \partial_2^P$ for every $K$-bundle 
$P \to X$. However, for $K = \PU(\cH)$, each 
element of $\Lambda^3(X,\Z)$ can be obtained as the Dixmier--Douady class  
of a $K$-bundle $P \to X$. 
\end{ex}

\begin{ex} We construct an example of a non-trivial aspherical 
obstruction class $\delta_1([P])$ for a central group extension 
with $\partial_2^{\hat K} \not=0$ and a space with non-trivial $\pi_3$. 

In view of Remark~\ref{rem:7.19}, 
it suffices to construct a compact topological space $X$ such 
that $H^3(X,\Z)\cong \Z$ and $\pi_3(X)$ is a torsion group because 
the latter assumption implies that 
$\partial_3^P \: \pi_3(X) \to \pi_2(K) \cong \Z$ 
vanishes for $K = \PU(\cH)$ and $\hat K = \U(\cH)$. 

We start with 
$A:=\bS^1\times\bS^2$, considered as the $3$-skeleton of a CW-complex. 
We attach to $A$ a $4$-cell $D^4$ via the attaching map
$f:\bS^3\to\bS^2,$ obtained 
by the composition of a map $\bS^3\to\bS^3$ of degree $p>1$ and the 
Hopf map $\bS^3\to\bS^2$ (which is the generator of $\pi_3(\bS^2)\cong\Z$).
We put
$$X:=\bS^1\times(\bS^2\cup_f D^4).$$

The space $X$ is then the $4$-cellular extension of $A$ in the sense of 
\cite[p.~48]{Whi}, and is compact by \cite[(1.1), p.~48]{Whi}.
The third homotopy group 
of $X$ can be computed using \cite[Thm (1.1) p.~212]{Whi}, 
which implies that the relative homotopy group $\pi_4(X,A)$ is a free 
$\Z$-module generated by the image of $f$, i.e. $\pi_4(X,A)\cong p\Z$. The 
long exact sequence of the pair now yields 
$$\pi_4(X,A)=p\Z\to\pi_3(A)=\Z\to\pi_3(X)\to\pi_3(X,A)=\{0\}$$ 
(\cite[Thm.~(2.4), p.~162]{Whi}). It follows that 
$\pi_3(X)\cong \Z/p\Z$ as claimed. 

The third cohomology of $X$ is computed in an indirect way. It is clear that
$H^*(A,\Z)$ is isomorphic to an exterior algebra in one generator 
$\theta$ of degree $1$ and one generator $\omega$ of degree $2$:
\[ H^*(A,\Z) = H^*(\bS^1\times\bS^2,\Z)\cong 
H^*(\bS^1,\Z)\otimes H^*(\bS^2,\Z)\cong 
\Lambda^*[\theta,\omega].\] 
Therefore 
$H^3(A,\Z)=\Z \theta \wedge \omega \cong \Z$. 
We only have to show that the attachment of the $4$-cell 
does not change $H^3$.    

The cohomology of a CW complex may be computed by cellular cohomology, see
\cite[p.~297]{Br97a}. In order to compute $H^3(X,\Z)$, we must add a 
$4$-cell in the complex $C^*(A,\Z)$ which computes the cellular cohomology
of $A$. But the coboundary operator for cellular cohomology 
is zero on the degree $4$-generator,
because the $4$-cell is attached to no $3$-cell (\cite[p.~297]{Br97a}). 
Therefore the attachment of 
the $4$-cell does not change $H^3$, which results in 
$H^3(X,\Z)\cong \Z$.  
\end{ex}

\section{Concluding remarks and problems} 

Proposition~\ref{prop:7.12} above gives a complete description of the 
aspherical obstruction classes for central extensions of $K$ 
by discrete groups. For central extensions by $K(\Gamma,1)$-groups 
our results are less complete. 
Theorem~\ref{thm:7.11} provides a solution for flat central 
extensions and bundles $P$ for which $\partial_2^P = 0$. 
The main open points concern the identification of aspherical obstructions 
on $1$-connected spaces $X$ with the corresponding central extensions 
of $\pi_2(X)$. Below we describe some concrete open problems concerning 
the concrete identification of obstruction classes.

\begin{prob}\mlabel{prob} It would be instructive to have more explicit 
  information on examples where $P = G$ is a Lie group, 
$K \subeq G$ is a closed subgroup and 
$X = G/K$ is a homogeneous space. Then we have for 
$Z = (\fz/\Gamma) \times D$ a natural map 
\[ \delta_1 \: \Ext_c(K,Z) \to 
H^3(G/K,\Gamma) \oplus H^2(G/K,D).\] 

(a) For the special case where $G$ is finite dimensional, 
the group $\pi_2(K)$ vanishes, so that the identity component 
of $Z$ leads to no non-trivial obstruction classes 
(Proposition~\ref{prop:0.1}). This leaves us with 
the case where $Z = D$ is discrete. If $K$ is connected, we thus obtain 
a homomorphism
\[ \delta_1 \: \Ext_c(K,Z) \cong \Hom(\pi_1(K),D) \to H^2(G/K,D).\] 
As $\pi_2(G)$ and $\pi_2(K)$ are trivial, we have an exact sequence 
\[ \{0\} \to \pi_2(G/K) \to \pi_1(K) \to \pi_1(G) \to \pi_1(G/K) \to 
\pi_0(K) \to \cdots. \]
If $G$ is $1$-connected and $K$ is connected, it follows that 
$G/K$ is $1$-connected with $H_2(G/K) \cong \pi_2(G/K) \cong \pi_1(K)$. 
This implies that $\delta_1$ is an isomorphism. 

(b) The case where $K \trile G$ is normal
is closely related to 
crossed modules (in the case that $K$ is closed and $G\to G/K$ admits a local continuous section) 
because for any lift of the $G$-action on 
$K$ to a central $Z$-extension $\hat K$, the homomorphism 
$p \: \hat K \to G$ defines a crossed module of groups whose characteristic 
class is an element of the third locally continuous
group cohomology 
$H^3_c(G/K,Z)$ (or locally smooth in the case of Lie groups). 
For an arbitrary topological group $Q$ there is a homomorphism
\begin{equation*}
\tau\: H^{3}_{c}(Q,Z)\to \check{H}^{2}(Q,\uline{Z}),
\end{equation*}
given as follows. If $f\: Q\times Q\times Q\to Z$ represents a class in $H^{3}_{c}(Q,Z)$,
then it is a group cocycle and continuous on $U\times U\times U$ for some open $U\subseteq Q$  with $e\in U$.
If $V\subseteq U$ is open, contains $e$ and satisfies 
$V^{2}\subseteq U$, then we have a \v{C}ech cocycle $\tau(f)_{g,h,k}\: g\cdot V\cap h\cdot V\cap k\cdot V\to Z$,
\begin{equation*}
 \tau(f)_{g,h,k}(x)=f(h,h^{-1}k,k^{-1}x)-f(g,g^{-1}k,k^{-1}x)+f(g,g^{-1}h,h^{-1}x),
\end{equation*}
for the open cover $(g\cdot V)_{g\in Q}$ of $Q$ (cf.\ \cite[Lemma 2.2]{Woc}).
Note that the cocycle identity for $f$ tells us that
\begin{equation*}
 \tau(f)_{g,h,k}(x)= f(g,g^{-1}h,h^{-1}k)+f(g^{-1}h,h^{-1}k,k^{-1}x)
\end{equation*}
depends continuously on $x$, since $g^{-1}h$, $h^{-1}k$ and $k^{-1}x$ are in 
$U$ if \break 
$g\cdot V\cap h\cdot V\cap k\cdot V\neq \emptyset$. Moreover, the class of
$\tau(f)$ clearly only depends on the class of $f$. 

If now $[c]\in H^{3}_{c}(G/K,Z)$ is the characteristic class
of the crossed module $\hat{K}\to G$ (\cite[Lemma 3.6]{Ne07}) 
and $G$ is considered as a $K$-principal bundle over $G/K$, then 
we expect that
\begin{equation*}
 \tau([c])=\pm \delta_{1}([G]), 
\end{equation*}
where $\hat{K}\to K$ is the central extension defining the crossed module.
\end{prob}



\begin{prob} (a) Consider $K$-bundles $q \: P \to X$ for which 
$\partial_2^P$ and $\partial_3^P$ vanish, so that 
$\pi_2(P)$ is an abelian extension of $\pi_2(X)$ by $\pi_1(K)$. 
Composing with $\partial_2^{\hat K}$, we obtain a class 
\[ (\partial_2^{\hat K})_*[\pi_2(P)] \in \Ext(\pi_2(X),\Gamma).\] 
If $X$ is $1$-connected, then  $\Lambda^3(X,\Gamma) \cong 
\Ext(\pi_2(X),\Gamma)$ (Remark~\ref{rem:6.1}(b)). Is it true that 
$(\partial_2^{\hat K})_*[\pi_2(P)]$ coincides with 
$\obs_P([\hat K])$? 

(b) If $Z$ is a $K(\Gamma,1)$-group and 
$\partial_2^{\hat K} \: \pi_2(K) \to \pi_1(Z)$ vanishes, 
then $\pi_1(\hat K)$ is an abelian extension of 
$\pi_1(K)$ by $\Gamma$ and for any $K$-bundle $P \to X$ 
we obtain a class
\[ (\partial_2^P)^*[\pi_1(\hat K)] \in \Ext(\pi_2(X),\Gamma).\] 
Suppose that $X$ is $1$-connected, which implies $\Lambda^3(X,\Gamma) \cong 
\Ext(\pi_2(X),\Gamma)$. Is it true that 
$(\partial_2^P)^*[\pi_1(\hat K)] = \obs_P([\hat K])$? 
\end{prob} 

\begin{prob} If $X$ is connected and $q_X \: \tilde X \to X$ is a universal covering, then 
$\pi_2(\tilde X) \cong \pi_2(X)$, and  we obtain 
a homomorphism 
\[ q_X^* \: \Lambda^3(X,\Gamma) \to \Lambda^3(\tilde X,\Gamma) \cong 
\Ext(\pi_2(X), \Gamma).\] 
It is an interesting problem to give an explicit description of the 
kernel of $q_X^*$ in $\Lambda^3(X,\Gamma)$. 
\end{prob}

\begin{prob} (Aspherical spaces)
Suppose that $X$ is aspherical, i.e., that $\tilde X$ is contractible. 
Then $X$ is a $K(\pi_1(X),1)$-space, so that 
\[ H^k(X,\Gamma) \cong H^k_{\rm grp}(\pi_1(X),\Gamma) \quad \mbox{ for } \quad 
k \in \N_0. \] 
Since all homotopy groups $\pi_k(X)$ vanish for $k \geq 2$, 
the corresponding connecting maps $\partial_k^P$ are trivial for any 
$K$-bundle $P$ over $X$. Therefore all obstruction classes are aspherical 
by Theorem~\ref{thm:main}. 

Write $q_X \: \tilde X \to X$ for the universal covering of $X$. 
Then, for every $K$-bundle $q \: P \to X$, the pullback 
$q_X^*P \to \tilde X$ is trivial because $\tilde X$ is contractible. 
As $q_X^*P \sim \tilde X \times K$ as $K$-bundles, the bundle $P$ can 
be written as a quotient 
\[ P \cong (\tilde X \times K)/\pi_1(X), \] 
where $\pi_1(X)$ acts on the trivial $K$-bundle $\tilde X \times K$ 
by bundle automorphism. 

Any lift of the $\pi_1(X)$-right action on $\tilde X$ to an action 
by bundle automorphisms on $\tilde X \times K$ is of the form 
\[ (x,k).d = (xd, f_d^{-1}(xd)k),\] 
where $f_d \in C(\tilde X, K)$ and the map 
\[ f \: \pi_1(X) \to C(\tilde X, K), \quad d \mapsto f_d  \] 
is a $1$-cocycle, i.e., 
\begin{equation}
  \label{eq:coc-cond}
f_{d_1 d_2} = f_{d_1} \cdot (d_1.f_{d_2}).
\end{equation}
We write $P^f := (\tilde X \times K)/\pi_1(X)$ for the quotient of 
$\tilde X \times K$ defined by the action given by the cocycle 
$f \: \pi_1(X) \to C(\tilde X,K)$ and  obtain a surjective map 
\[ H^1_{\rm grp}(\pi_1(X), C(\tilde X,K)) \to \Bun(X,K), \quad 
[f] \mapsto [P^f].\] 

Let $\hat K$ be a topological $Z$-extension of $K$. 
Then we have a natural map 
$C(\tilde X,\hat K) \to C(\tilde X,K)$ 
which is surjective because, for every continuous function 
$f \: \tilde X \to K$, the $Z$-bundle $f^*\hat K \to \tilde X$ is 
trivial. We thus have a central extension 
\[ \1 \to C(\tilde X,Z) \to C(\tilde X,\hat K) \to C(\tilde X,K) \to \1 \] 
which leads to long exact sequence in group cohomology. In particular, 
if 
\[ \delta_1 \: H^1_{\rm grp}(\pi_1(X), C(\tilde X,K)) 
\to H^2_{\rm grp}(\pi_1(X),C(\tilde X,Z)).\] 
is the connecting map, then $\delta_1([f])$ vanishes 
if and only if the $K$-bundle $P^f$ lifts to a $\hat K$-bundle. 
Next we consider the short exact sequence 
\[ \1 \to C(\tilde X,\Gamma) \cong \Gamma \to 
C(\tilde X,\tilde Z) \to C(\tilde X,Z) \to \1 \]
which in turn leads to a connecting map 
\[ \delta_2 \: H^2_{\rm grp}(\pi_1(X), C(\tilde X,Z)) 
\to H^3_{\rm grp}(\pi_1(X),\Gamma) \cong H^3(X,\Gamma).\] 
Is it true that 
\[ \delta_1([P^f]) = \pm \delta_2(\delta_1([f]))? \]
\end{prob}

\end{document}